\documentclass[11pt,letterpaper,reqno]{amsart}
\usepackage{amsmath,amsthm,amsfonts,amssymb,mathtools,algpseudocode}
\usepackage{comment,bm,color,mathdots,bookmark,tikz}
\usepackage{caption,subcaption}

\hypersetup{pdfstartview={FitH}}

\newtheorem{Theorem}{{\bf Theorem}}[section]

\newtheorem{Algorithm}{{\bf Algorithm}}

\newtheorem{Lemma}[Theorem]{{\bf Lemma}}
\newtheorem{Remark}[Theorem]{{\bf Remark}}

\newcommand{\diag}{\textnormal{diag}}
\newcommand{\dd}{\emph{diag}}

\newcommand{\grad}{\textnormal{grad}}

\newcommand{\tr}{\textnormal{trace}}
\newcommand{\im}{\textnormal{Im}}
\newcommand{\re}{\textnormal{Re}}
\newcommand{\N}{\mathcal{N}}
\newcommand{\U}{\mathcal{U}}
\newcommand{\Ham}{\mathcal{H}}
\newcommand{\sHam}{\mathcal{W}}
\newcommand{\perH}{\mathcal{M}}
\newcommand{\sperH}{\mathcal{K}}
\newcommand{\Symp}{\mathcal{S}p}
\newcommand{\Perp}{\mathcal{P}p}
\numberwithin{equation}{section}

\begin{document}

\title[Finding the closest normal structured matrix]{Finding the closest normal structured matrix}

%\author{One author\thanks{some info} \and Another author\thanks{more info}}
\author{Erna Begovi\'{c}~Kova\v{c}}\thanks{\textsc{Erna Begovi\'{c}~Kova\v{c},
Faculty of Chemical Engineering and Technology, University of Zagreb, Maruli\'{c}ev trg 19, 10000 Zagreb, Croatia}.
\texttt{ebegovic@fkit.hr}}
%\thanks{}
\date{\today}

\renewcommand{\subjclassname}{\textup{2020} Mathematics Subject Classification}

\begin{abstract}
Given a structured matrix $A$ we study the problem of finding the closest normal matrix with the same structure.
The structures of our interest are: Hamiltonian, skew-Hamiltonian, per-Hermitian, and perskew-Hermitian.
We develop a structure-preserving Jacobi-type algorithm for finding the closest normal structured matrix and show that such algorithm converges to a stationary point of the objective function.
\end{abstract}

\maketitle

\section{Introduction}

The problem of finding the closest normal matrix $X$ to any unstructured matrix $A\in\mathbb{C}^{n\times n}$ in the Frobenius norm
\begin{equation}\label{minimization}
\min_{X \in \N}\|X-A\|_F^2,
\end{equation}
where $\N$ stands for the set of normal matrices,
was an open question for a long time.
It was solved independently by Gabriel~\cite{Gabriel79,Gabriel87} and Ruhe~\cite{Ruhe87}. A nice summary of important
findings is given by Higham in~\cite{Higham89}. In this paper we are interested in the structure-preserving version of problem~\eqref{minimization}. That is, given a structure $\mathcal{S}$ and matrix $A\in\mathcal{S}$, we are looking for
\begin{equation}\label{minimization-S}
\min_{X\in\N\cap\mathcal{S}}\|X-A\|_F.
\end{equation}

The following theorem from~\cite{CauseyPhD} states the solution of~\eqref{minimization} using a maximization problem formulation.
See~\cite[Theorem 5.2]{Higham89} for a full set of references. Notation $\U$ stands for the set of unitary matrices.

\begin{Theorem}\label{tm:maximization-normal}
Let $A\in\mathbb{C}^{n\times n}$ and let $X=UDU^H$, where $U\in\U$ and $D\in\mathbb{C}^{n\times n}$ is diagonal. Then $X$ is a nearest normal matrix to $A$ in the Frobenius norm if and only if
\begin{itemize}
\item[(a)] $\|{\dd}(U^HAU)\|_F=\displaystyle{\max_{Q\in\U}}\|{\dd}(Q^HAQ)\|_F$, and
\item[(b)] $ D={\dd}(U^HAU)$.
\end{itemize}
\end{Theorem}

Thus, the problem of finding the closest normal matrix to $A\in\mathbb{C}^{n\times n}$ can be transformed into a problem of finding a unitary similarity transformation $Q$ which makes the sum of squares of the diagonal elements of $Q^HAQ$ as large as possible. Instead of solving the minimization problem~\eqref{minimization}, one can address the dual maximization problem
\begin{equation}\label{maximization}
\max_{Q\in\U}\|\diag(Q^HAQ)\|_F^2.
\end{equation}

It is well known that if $A$ is normal, then it can be unitarily diagonalizable. Since we focus on matrices that are not normal, the goal is to make the matrix $Q^HAQ$ ``as diagonal as possible''. Then, the closest normal matrix to $A$ is obtained as $X =Q\diag(Q^HAQ)Q^H$.

We consider four classes of matrices:
\begin{itemize}
\item Hamiltonian $\Ham=\{A\in\mathbb{C}^{2n\times2n} \ | \ (JA)^H=JA\}$,
\item skew-Hamiltonian $\sHam=\{A\in\mathbb{C}^{2n\times2n} \ | \ (JA)^H=-JA\}$,
\item per-Hermitian $\perH=\{A\in\mathbb{C}^{m\times m} \ | \ (FA)^H=FA\}$,
\item perskew-Hermitian $\sperH=\{A\in\mathbb{C}^{m\times m} \ | \ (FA)^H=-FA\}$,
\end{itemize}
where
\begin{equation}\label{JF}
J=J_{2n}=\left[
      \begin{array}{cc}
        0 & I_n \\
        -I_n & 0 \\
      \end{array}
    \right]\in\mathbb{R}^{2n\times2n}, \quad
F=F_m=\left[
\begin{array}{ccc}
&&1\\
&\iddots&\\
1&&
  \end{array}
\right]\in\mathbb{R}^{m\times m}.
\end{equation}
A unitary similarity transformation $Z^HAZ$, $Z\in\U$ is, in  general, not structure-preserving. Therefore, in order to get $Z^HAZ\in\mathcal{S}$ for $A\in\mathcal{S}$, matrix $Z$ needs to have an additional structure.
Transformations that keep the structure of the sets $\Ham$ and $\sHam$ are symplectic transformations
$$\Symp=\{Z\in\mathbb{C}^{2n\times2n} \ | \ Z^HJZ=J\},$$
and transformations that keep the structure of the sets $\perH$ and $\sperH$ are perplectic transformations
$$\Perp=\{Z\in\mathbb{C}^{m\times m} \ | \ Z^HFZ=F\}.$$

Both groups of symplectic and perplectic matrices form manifolds. Hamiltonian matrices form the tangent subspace on the manifold of symplectic matrices at the identity. It is easy to check this. For symplectic matrices $M$ we have $h(M):=M^HJM-J=0$. Tangent space at the identity is the set of matrices $\{A \ | \ Dh(I)A=0\}$. Using linear approximation we get
$$h(I+A)=h(I)+Dh(I)A+O(\|A\|^2),$$
and since $h(I)=0$,
$$Dh(I)A=(I+A)^HJ(I+A)-J=A^HJ+JA+O(\|A\|^2)=0.$$
Matrices that satisfy the equation $A^HJ+JA=0$ are indeed Hamiltonian matrices. Orthogonal space at the identity is orthogonal complement of the set of Hamiltonian matrices, which is the set of skew-Hamiltonain matrices. In the same way one can check that perskew-Hermitian matrices form the tangent subspace on the manifold of perplectic matrices at the identity, and per-Hermitian matrices form its orthogonal subspace.
Transformations from a manifold preserve the structure of the matrices from the corresponding tangent or orthogonal subspace.

In the algebraic setting, one can look at the symplectic and perplectic groups as Lie groups. Hamiltonian and skew-Hamiltonian matrices are Lie algebra and Jordan algebra of the symplectic group, respectively, while per-Hermitian and perskew-Hermitian matrices are Jordan algebra and Lie algebra of the perplectic group, respectively. Transformations from a Lie group preserve the structure of the corresponding Jordan or Lie algebra.

Both geometric and algebraic interpretation of the studied matrix structures are given in Table~\ref{table:geomalg}. One can also find more about these structures in the existing literature, e.g.,~\cite{MMT03,Trench04}.
\bigskip

\begin{table}[h!]
\centering
\begin{tabular}{|c|c|c|}
  \hline
  manifold & tangent subspace at $I$ & orthogonal subspace at $I$ \\
  \hline
  \hline
  symplectic & Hamiltonian & skew-Hamiltonian \\
  perplectic & perskew-Hermitian & per-Hermitian \\
  \hline
  \hline
  Lie group & Lie algebra & Jordan algebra \\
  \hline
\end{tabular}
\caption{Geometric and algebraic setting for the structured matrices}
\label{table:geomalg}
\end{table}

In Section~\ref{sec:theorems} we give structured analogues of Theorem~\ref{tm:maximization-normal} and formulate the corresponding versions of minimization problem~\eqref{minimization-S}. Then in Section~\ref{sec:algorithm} we develop the Jacobi-type algorithm for solving the minimization problems defined in Section~\ref{sec:theorems} and prove its convergence in Section~\ref{sec:cvgproof}. Finally, in Section~\ref{sec:num} we present some numerical results.

\section{Structured analogues of Theorem~\ref{tm:maximization-normal}}\label{sec:theorems}

We study minimization problem~\eqref{minimization-S}. Theorem~\ref{tm:maximization-normal} suggests to find a unitary matrix $U$ that maximizes $\|\diag(U^HAU)\|_F$. Let us explore how that approach can be used with the structure-preserving constrain.

\subsection{Hamiltonian and skew-Hamiltonian matrices}

A Hamiltonian matrix  $H\in\mathbb{C}^{2n\times2n}$ can be written as a $2\times2$ block matrix
\begin{equation}\label{eq:Hamsetup}
H=\left[
      \begin{array}{cc}
        H_{11} & H_{12} \\
        H_{21} & -H_{11}^H \\
      \end{array}
    \right], \quad \textnormal{where} \ H_{12}^H=H_{12}, \ H_{21}^H=H_{21}, \quad H_{11},H_{12},H_{21}\in\mathbb{C}^{n\times n}.
\end{equation}
Moreover, a skew-Hamiltonian matrix $W\in\mathbb{C}^{2n\times2n}$ can be written as a $2\times2$ block matrix
\begin{equation}\label{eq:sHamsetup}
W=\left[
      \begin{array}{cc}
        W_{11} & W_{12} \\
        W_{21} & W_{11}^H \\
      \end{array}
    \right], \quad \text{where} \ W_{12}^H=-W_{12}, \ W_{21}^H=-W_{21}, \quad W_{11},W_{12},W_{21}\in\mathbb{C}^{n\times n}.
\end{equation}
It follows from~\eqref{eq:Hamsetup} and ~\eqref{eq:sHamsetup}, respectively, that any diagonal Hamiltonian matrix has the form
$$D_H=\left[
       \begin{array}{cc}
         D & 0 \\
         0 & -D^H \\
       \end{array}
     \right],$$
while any skew-Hamiltonian diagonal matrix has the form
$$D_W=\left[
       \begin{array}{cc}
         D & 0 \\
         0 & D^H \\
       \end{array}
     \right],$$
where $D=\diag(d_1,\ldots,d_n)$.
Also, it is easy to check that for every skew-Hamiltonian matrix $W\in\sHam$ there is a Hamiltonian matrix $H\in\Ham$
(and for every $H\in\Ham$ there is $W\in\sHam$) such that
$$W=\imath H.$$
Therefore, all results obtained for Hamiltonian matrices will imply analogue results for skew-Hamiltonian matrices.

In order to obtain a result analogue to that in Theorem~\ref{tm:maximization-normal}, we use the Schur decomposition for Hamiltonian matrices given in~\cite{PvL81}.

\begin{Theorem}[\cite{PvL81}]\label{tm:HamSchur}
If $H\in\mathbb{C}^{2n\times2n}$ is a Hamiltonian matrix whose eigenvalues have nonzero real parts, then there exists a unitary
$$U=\left[
      \begin{array}{cc}
        U_{11} & U_{12} \\
        -U_{12} & U_{11} \\
      \end{array}
    \right], \qquad U_{11}, U_{12}\in\mathbb{C}^{n\times n},
$$
such that
\begin{equation}\label{eq:HamSchur}
U^HHU=\left[
          \begin{array}{cc}
            T & M \\
            0 & -T^H \\
          \end{array}
        \right], \qquad T, M\in\mathbb{C}^{n\times n},
\end{equation}
where $T$ is upper triangular and $M^H=M$.
\end{Theorem}

The following lemma is a special case of Theorem~\ref{tm:HamSchur} for Hamiltonian and skew-Hamiltonian normal matrices.

\begin{Lemma}\label{lemma:Schur-nH}
\begin{itemize}
\item[(i)] If $H\in\mathbb{C}^{2n\times2n}$ is a normal Hamiltonian matrix whose eigenvalues have nonzero real parts, then there exists a unitary symplectic $U\in\mathbb{C}^{2n\times2n}$ and diagonal $D\in\mathbb{C}^{n\times n}$ such that
\begin{equation}\label{eq:Schur-nH}
H=U\left[
       \begin{array}{cc}
         D & 0 \\
         0 & -D^H \\
       \end{array}
     \right]U^H.
\end{equation}
\item[(ii)] If $W\in\mathbb{C}^{2n\times2n}$ is a normal skew-Hamiltonian matrix whose eigenvalues have nonzero imaginary parts, then there exists a unitary symplectic $U\in\mathbb{C}^{2n\times2n}$ and diagonal $D\in\mathbb{C}^{n\times n}$ such that
\begin{equation}\label{eq:Schur-nsH}
W=U\left[
       \begin{array}{cc}
         D & 0 \\
         0 & D^H \\
       \end{array}
     \right]U^H.
\end{equation}
\end{itemize}
\end{Lemma}

\begin{proof}
\begin{itemize}
\item[(i)] The Schur decomposition of $H\in\Ham$ is as in relation~\eqref{eq:HamSchur}. Matrices $H$ and $U^HHU$ are normal. For permutation
$P=\left[
      \begin{array}{cc}
        I_n &  \\
         & F_n \\
      \end{array}
    \right]$,
where $F_n$ is as in~\eqref{JF}, matrix
\begin{equation}\label{Schur-nH1}
P^T(U^HHU)P=P^T\left[
          \begin{array}{cc}
            T & M \\
            0 & -T^H \\
          \end{array}
        \right]P
\end{equation}
is normal and triangular.
Normal triangular matrix must be diagonal. Diagonal elements of $P^T(U^HHU)P$ are the same as of $U^HHU$.
Therefore, in~\eqref{Schur-nH1} we conclude that $M=0$, $T=D$ is diagonal, and
$$U^HHU=\left[
          \begin{array}{cc}
            D & 0 \\
            0 & -D^H \\
          \end{array}
        \right].$$
This gives relation~\eqref{eq:Schur-nH}.

It is easy to check that matrix $U$ is indeed symplectic. We have
\begin{align*}
JU & =\left[
      \begin{array}{cc}
        0 & I \\
        -I & 0 \\
      \end{array}
    \right]\left[
      \begin{array}{cc}
        U_{11} & U_{12} \\
        -U_{12} & U_{11} \\
      \end{array}
    \right]=\left[
      \begin{array}{cc}
        -U_{12} & U_{11} \\
        -U_{11} & -U_{12} \\
      \end{array}
    \right] \\
& =\left[
      \begin{array}{cc}
        U_{11} & U_{12} \\
        -U_{12} & U_{11} \\
      \end{array}
    \right]\left[
      \begin{array}{cc}
        0 & I \\
        -I & 0 \\
      \end{array}
    \right]=UJ.
\end{align*}
Now, since $U$ is unitary, it follows that $U^HJU=J$.

\item[(ii)] Let $W\in\sHam$. Then $W=\imath H$ for some $H\in\Ham$. If $W$ is normal, then $H$ is also normal. If the eigenvalues of $W$ have nonzero imaginary parts, then eigenvalues of $H$ have nonzero real parts. Hence, we can apply first assertion of this lemma on $H$. This gives
$$H=U\left[
          \begin{array}{cc}
            \tilde{D} & 0 \\
            0 & -\tilde{D}^H \\
          \end{array}
        \right]U^H.$$
For $D=\imath\tilde{D}$ it follows
$$W=\imath H=U\left[
       \begin{array}{cc}
         \imath\tilde{D} & 0 \\
         0 & -\imath\tilde{D}^H \\
       \end{array}
     \right]U^H=
     U\left[
       \begin{array}{cc}
         D & 0 \\
         0 & D^H \\
       \end{array}
     \right]U^H.$$
\end{itemize}
\end{proof}

Using the decompositions from Lemma~\ref{lemma:Schur-nH} we will prove Theorems~\ref{tm:maximizationH} and~\ref{tm:maximizationsH} which are structured analogues to Theorem~\ref{tm:maximization-normal} for structures $\Ham$ and $\sHam$, respectively.
Before that, we need one more auxiliary result.

\begin{Lemma}\label{lema:pitagora}
For a general matrix $M$ we have
\begin{equation}\label{pitagora}
\|M\|^2=\|M-\diag(M)\|^2+\|\diag(M)\|^2.
\end{equation}
\end{Lemma}

\begin{proof}
Let $M$ be an arbitrary matrix. Its orthogonal projection to the subspace of diagonal matrices is $\diag(M)$. On the other hand, null-matrix $\mathbf{0}$ also belongs to the subspace of the diagonal matrices. Hence, matrices $M$, $\diag(M)$ and $\mathbf{0}$ are vertices of a right-angled triangle with legs $\|\diag(M)\|$ and $\|M-\diag(M)\|$ and the hypothenuse $\|M\|$.

\begin{center}
\begin{tikzpicture}
\draw [-] (-4,0) -- (8,0);
\node [below] at (7,0) {diagonal matrices};
\draw [ultra thick,blue] (0,0) -- (4,0) -- (0,2) -- (0,0);
\draw[fill,blue] (0,0) circle [radius=0.05];
\draw[fill,blue] (4,0) circle [radius=0.05];
\draw[fill,blue] (0,2) circle [radius=0.05];
\node [below,blue] at (0,0) {$\diag(M)$};
\node [below,blue] at (4,0) {$\mathbf{0}$};
\node [above,blue] at (0,2) {$M$};
\node [left] at (0,0.75) {$\|M-\diag(M)\|$};
\node [right] at (1.75,1.25) {$\|M\|$};
\node [below] at (2.5,0) {$\|\diag(M)\|$};
\end{tikzpicture}
\end{center}

Now, equation~\eqref{lema:pitagora} follows from the Pythagoras' theorem.
\end{proof}

\begin{Theorem}\label{tm:maximizationH}
Let $A\in\mathbb{C}^{2n\times2n}$ be a Hamiltonian matrix and let $X=ZDZ^H$, where $Z$ is symplectic unitary and $D$ is Hamiltonian diagonal. Then $X$ is a normal Hamiltonian matrix with no purely imaginary eigenvalues, closest to $A$ in the Frobenius norm, if and only if
\begin{itemize}
\item[(a)] $\displaystyle \|\emph{\diag}(Z^HAZ)\|_F=\max_{Q\in\U\cap\Symp}\|\emph{\diag}(Q^HAQ)\|_F$, and
\item[(b)] $\displaystyle D=\emph{\diag}(Z^HAZ)$.
\end{itemize}
\end{Theorem}

\begin{proof}
Let $A\in\Ham$. By $X_0$ denote the closest normal Hamiltonian matrix to $A$. If $A$ is already normal, the distance between $A$ and $X_0$ is zero. Otherwise,
\begin{equation}\label{closestnH}
\min_{X\in\N\cap\Ham}\|A-X\|_F=\|A-X_0\|_F.
\end{equation}

Let $X_0=ZD_0Z^H$ be the Schur decomposition of $X_0\in\N\cap\Ham$, like in~\eqref{eq:Schur-nH}, $D_0\in\Ham$.
Then
$$\|A-X_0\|_F=\|A-ZD_0Z^H\|_F=\|Z^HAZ-D_0\|_F,$$
and relation~\eqref{closestnH} is transformed into
$$\min_{D\in\Ham \ \text{diagonal}}\|Z^HAZ-D\|_F=\|Z^HAZ-D_0\|_F.$$
The closest diagonal matrix to $Z^HAZ$ is its orthogonal projection to the subspace of diagonal matrices, which is simply $\diag(Z^HAZ)$.
This gives $D_0=\diag(Z^HAZ)$ and implies assertion $(b)$.

To obtain $(a)$, take $D=\diag(Q^HAQ)\in\Ham$, $Q$ unitary symplectic. Matrix $Q(\diag(Q^HAQ))Q^H$ is normal and its distance from $A$ is at least $X_0$. Thus,
\begin{align}
\|A-Q(\diag(Q^HAQ))Q^H\|_F^2 & \geq\|A-X_0\|_F^2, \nonumber \\
\|Q^HAQ-\diag(Q^HAQ)\|_F^2 & \geq\|Z^HAZ-Z^HX_0Z\|_F^2. \label{pitagora1}
\end{align}
On the left-hand side of~\eqref{pitagora1} we use Lemma~\eqref{lema:pitagora} for $M=Q^HAQ$, while on the right-hand side we use the same lemma for $M=Z^HAZ$. We get
\begin{align*}
\|Q^HAQ\|_F^2-\|\diag(Q^HAQ)\|_F^2 & \geq\|Z^HAZ\|_F^2-\|Z^HX_0Z\|_F^2, \\
\|\diag(Q^HAQ)\|_F^2 & \leq\|Z^HX_0Z\|_F^2, \\
\max_{Q\in\U\cap\Symp}\|\diag(Q^HAQ)\|_F^2 & =\|Z^HX_0Z\|_F^2.
\end{align*}

Conversely, let $(a)$ and $(b)$ hold for $X_0\in\N\cap\Ham$. There exists a closest normal Hamiltonian matrix because both set of normal and set of Hamiltonian matrices are closed. Assume that $X_0$ is not the closest, that is
$$\min_{X\in\N\cap\Ham}\|A-X\|_F\neq\|A-X_0\|_F.$$
Then $\|A-X\|<\|A-X_0\|,$ for some $X\in\N\cap\Ham$.
Take $X_0=ZD_0Z^H$, $X=QDQ^H$ from the Schur decomposition~\eqref{eq:Schur-nH}. It follows from $(b)$ that $D_0=\diag(Z^HAZ)$ and $D=\diag(Q^HAQ)$.
Using the argument~\eqref{pitagora} again, we get
\begin{align*}
\|A-Q\diag(Q^HAQ)Q^H\|_F^2 & <\|A-Z\diag(Z^HAZ)Z^H\|_F^2, \\
\|Q^HAQ-\diag(Q^HAQ)\|_F^2 & <\|Z^HAZ-\diag(Z^HAZ)\|_F^2, \\
\|Q^HAQ\|_F^2-\|\diag(Q^HAQ)\|_F^2 & <\|Z^HAZ\|_F^2-\|\diag(Z^HAZ)\|_F^2, \\
\|\diag(Q^HAQ)\|_F^2 & >\|\diag(Z^HAZ)\|_F^2, \\
\max_{Q\in\U\cap\Symp}\|\diag(Q^HAQ)\|_F^2 & >\|\diag(Z^HAZ)\|_F^2,
\end{align*}
which is contradiction with $(a)$.
\end{proof}

\begin{Theorem}\label{tm:maximizationsH}
Let $A\in\mathbb{C}^{2n\times2n}$ be a skew-Hamiltonian matrix and let $X=ZDZ^H$, where $Z$ is symplectic unitary and $D$ is skew-Hamiltonian diagonal. Then $X$ is a normal skew-Hamiltonian matrix with no real eigenvalues, closest to $A$ in the Frobenius norm, if and only if
\begin{itemize}
\item[(a)] $\displaystyle \|\emph{\diag}(Z^HAZ)\|_F=\max_{Q\in\U\cap\Symp}\|\emph{\diag}(Q^HAQ)\|_F$, and
\item[(b)] $\displaystyle D=\emph{\diag}(Z^HAZ)$.
\end{itemize}
\end{Theorem}

\begin{proof}
The proof is the same as for the Hamiltonian case from Theorem~\ref{tm:maximizationH}, but instead of the Schur decomposition~\eqref{eq:Schur-nH} it uses~\eqref{eq:Schur-nsH}.
\end{proof}

\subsection{Per-Hermitian and perskew-Hermitian matrices}

A per-Hermitian matrix  $M\in\mathbb{C}^{2n\times2n}$ can be written as a $2\times2$ block matrix
$$M = \left[ \begin{array}{cc}
M_{11} & M_{12}\\
M_{21} & FM_{11}^HF
\end{array}\right], \quad \textnormal{where} \ (FM_{12})^H=FM_{12}, \ (FM_{21})^H=FM_{21},$$
$M_{11},M_{12},M_{21}\in\mathbb{C}^{n\times n}$. The elements of the antidiagonal of $M_{12}$ and $M_{21}$ have to be real.
A perskew-Hermitian matrix $K\in\mathbb{C}^{2n\times2n}$ can be written as a $2\times2$ block matrix
$$K = \left[ \begin{array}{cc}
K_{11} & K_{12}\\
K_{21} & -FK_{11}^HF
\end{array}\right], \quad \textnormal{where} \ (FK_{12})^H=-FK_{12}, \ (FK_{21})^H=-FK_{21},$$
$K_{11},K_{12},K_{21}\in\mathbb{C}^{n\times n}$. The elements of the antidiagonal of $K_{12}$ and $K_{21}$ have to be imaginary (or zero).
A diagonal per-Hermitian matrix $D_M$ has to be of the form
$$D_M=\left[\begin{array}{cc}
D & 0\\
0 & FD^HF
\end{array}\right],$$
while a diagonal perskew-Hermitian matrix $D_K$ is given by
$$D_K=\left[\begin{array}{cc}
D & 0\\
0 & -FD^HF
\end{array}\right],$$
where $D=\diag(d_1,\ldots,d_n)$.
Also, for every perskew-Hermitian matrix $K\in\sperH$ there is a per-hermitian matrix $M\in\perH$, and viceversa, such that
$$K=\imath M.$$

The next lemma gives Schur-like decomposition of per-Hermitian and perskew-Hermitian normal matrices.

\begin{Lemma}\label{lemma:Schur-nper}
\begin{itemize}
\item[(i)] If $A\in\mathbb{C}^{2n\times2n}$ is a normal per-Hermitian matrix whose eigenvalues have nonzero imaginary parts, then there exists a unitary perplectic $U\in\mathbb{C}^{2n\times2n}$ and diagonal $D\in\mathbb{C}^{n\times n}$ such that
\begin{equation}\label{eq:Schur-nper}
A=U\left[
       \begin{array}{cc}
         D & 0 \\
         0 & RD^HR \\
       \end{array}
     \right]U^H.
\end{equation}
\item[(ii)] If $A\in\mathbb{C}^{2n\times2n}$ is a normal perskew-Hermitian matrix whose eigenvalues have nonzero real parts, then there exists a unitary perplectic $U\in\mathbb{C}^{2n\times2n}$ and diagonal $D\in\mathbb{C}^{n\times n}$ such that
\begin{equation}\label{eq:Schur-npers}
A=U\left[
       \begin{array}{cc}
         D & 0 \\
         0 & -RD^HR \\
       \end{array}
     \right]U^H.
\end{equation}
\end{itemize}
\end{Lemma}

\begin{proof}
\begin{itemize}
\item[(i)]
First, notice that eigenvalues of per-Hermitian matrix $A$ come in complex conjugate pairs $(\lambda,\bar{\lambda})$ with $\lambda$ and $\bar{\lambda}$ having the same algebraic multiplicity. Let us verify this. If $\lambda\in\sigma(A)$, then $\bar{\lambda}\in\sigma(A^H)$. Since $\sigma(A^H)=\sigma(RA^HR)$ and $RA^HR=A$, we have $\bar{\lambda}\in\sigma(A)$.

Let $\{\lambda_1,\ldots,\lambda_p,\bar{\lambda}_1,\ldots,\bar{\lambda}_p\}$ be the eigenvalues of $A$ and let $v_1,\ldots,v_n$ be a complete set of orthogonal eigenvectors corresponding to $\lambda_1,\ldots,\lambda_p$.
Set $V=\left[
     \begin{array}{cccc}
       v_1 & v_2 & \cdots & v_n \\
     \end{array}
   \right]\in\mathbb{C}^{2n\times n}$.
If $v_i,v_j\in\mathbb{C}^{2n}$ are eigenvectors of $A$ for $\lambda_i$ and $\lambda_j$, respectively, then $v_i^HR_{2n}v_j\neq0$ only if $\lambda_i=\bar{\lambda}_j$. For $i\neq j$ we have $\lambda_i\neq\bar{\lambda}_j$ and
since all eigenvalues of $A$ have nonzero imaginary parts, we have
$\lambda\neq\bar{\lambda}$ for all $\lambda\in\sigma(A)$. This implies that
\begin{equation}\label{Schur-nper1}
V^HR_{2n}V=0.
\end{equation}

Define $U:=\left[
     \begin{array}{cc}
       V & R_{2n}VR_n \\
     \end{array}
   \right]\in\mathbb{C}^{2n\times 2n}.$
Using identity~\eqref{Schur-nper1} along with $V^HV=I_n$ and $RR=I$ it is easy to check that $U$ is unitary
\begin{align*}
U^HU & =\left[
       \begin{array}{c}
         V^H \\
         R_nV^HR_{2n} \\
       \end{array}
     \right]\left[
     \begin{array}{cc}
       V & R_{2n}VR_n \\
     \end{array}
   \right] \\
& =\left[
  \begin{array}{cc}
    V^HV & V^HR_{2n}VR_n \\
    R_nV^HR_{2n}V & R_nV^HR_{2n}R_{2n}VR_n \\
  \end{array}
\right]=\left[
     \begin{array}{cc}
       I_n & 0 \\
       0 & I_n \\
     \end{array}
   \right]=I_{2n},
\end{align*}
and perplectic
\begin{align*}
U^HR_{2n}U & =\left[
       \begin{array}{c}
         V^H \\
         R_nV^HR_{2n} \\
       \end{array}
     \right]R_{2n}\left[
     \begin{array}{cc}
       V & R_{2n}VR_n \\
     \end{array}
   \right] \\
& =\left[
  \begin{array}{cc}
    V^HR_{2n}V & V^HR_{2n}R_{2n}VR_n \\
    R_nV^HR_{2n}R_{2n}V & R_nV^HR_{2n}R_{2n}R_{2n}VR_n \\
  \end{array}
\right]=\left[
     \begin{array}{cc}
       0 & R_n \\
       R_n & 0 \\
     \end{array}
   \right]=R_{2n}.
\end{align*}

Then, using the fact that $\text{span}(V)$ is invariant subspace for $A$, that $AV=VB$ for some $B$, and $RAR=A^H$, we have
\begin{align*}
U^HAU & =\left[
     \begin{array}{cc}
       V^HAV & V^HAR_{2n}VR_n \\
       R_nV^HR_{2n}AV & R_nV^HR_{2n}AR_{2n}VR_n \\
     \end{array}
   \right] \\
& =\left[
     \begin{array}{cc}
       V^HAV & V^HAR_{2n}VR_n \\
       R_nV^HR_{2n}VB & R_nV^HA^HVR_n \\
     \end{array}
   \right] = \left[
     \begin{array}{cc}
       D & V^HAR_{2n}VR_n \\
       0 & R_nD^HR_n \\
     \end{array}
   \right],
\end{align*}
for $D=\text{diag}(\lambda_1,\ldots,\lambda_p)$.
Since $A$ is assumed to be normal, $U^HAU$ is upper triangular and normal. Therefore, it is diagonal, that is $V^HAR_{2n}VR_n=0$,
which gives decomposition~\eqref{eq:Schur-nper}.

\item[(ii)] If $A$ is perskew-Hermitian, its eigenvalues come in pairs $(\lambda,-\bar{\lambda})$ and the assumption of nonzero real parts assures that $\lambda\neq-\bar{\lambda}$.
Further on, the proof follows the same reasoning as above.
\end{itemize}
\end{proof}

We use decompositions from Lemma~\ref{lemma:Schur-nper} to get the results analogue to those in Theorems~\ref{tm:maximizationH} and~\ref{tm:maximizationsH}.

\begin{Theorem}\label{tm:maximizationper}
Let $A\in\mathbb{C}^{2n\times2n}$ be a per-Hermitian matrix and let $X=ZDZ^H$, where $Z$ is perplectic unitary and $D$ is per-Hermitian diagonal. Then $X$ is a normal per-Hermitian matrix with no real eigenvalues, closest to $A$ in the Frobenius norm, if and only if
\begin{itemize}
\item[(a)] $\displaystyle \|\diag(Z^HAZ)\|_F=\max_{Q\in\U\cap\Perp}\|\diag(Q^HAQ)\|_F$, and
\item[(b)] $\displaystyle D=\diag(Z^HAZ)$.
\end{itemize}
\end{Theorem}

\begin{proof}
The proof follows the lines of the proof of Theorem~\ref{tm:maximizationH}. Instead of symplectic we have perplectic matrices and instead of decomposition~\eqref{eq:Schur-nH} we use~\eqref{eq:Schur-nper}.
\end{proof}

\begin{Theorem}\label{tm:maximizationpers}
Let $A\in\mathbb{C}^{2n\times2n}$ be a perskew-Hermitian matrix and let $X=ZDZ^H$, where $Z$ is perplectic unitary and $D$ is perskew-Hermitian diagonal. Then $X$ is a normal perskew-Hermitian matrix with no purely imaginary eigenvalues, closest to $A$ in the Frobenius norm, if and only if
\begin{itemize}
\item[(a)] $\displaystyle \|\diag(Z^HAZ)\|_F=\max_{Q\in\U\cap\Perp}\|\diag(Q^HAQ)\|_F$, and
\item[(b)] $\displaystyle D=\diag(Z^HAZ)$.
\end{itemize}
\end{Theorem}

\begin{proof}
The proof follows the lines of the proof of Theorem~\ref{tm:maximizationH}. Instead of symplectic we have perplectic matrices and instead of decomposition~\eqref{eq:Schur-nH} we use~\eqref{eq:Schur-npers}.
\end{proof}

\section{Jacobi-type algorithm for finding the closest normal matrix with a given structure}\label{sec:algorithm}

Based on the results from Section~\ref{sec:theorems} we can formulate the structured analogues of the maximization problem~\eqref{maximization}.
Assuming that $A$ is Hamiltonian or skew-Hamiltonian, it follows from Theorems~\ref{tm:maximizationH} and~\ref{tm:maximizationsH}, respectively, that a dual maximization formulation of the minimization problem~\eqref{minimization-S} is
\begin{equation}\label{maximizationS}
\max_{Z\in\U\cap\Symp}\|\diag(Z^HAZ)\|_F^2,
\end{equation}
while in per-Hermitian or perskew-Hermitian case Theorems~\ref{tm:maximizationper} and~\ref{tm:maximizationpers} imply the form
\begin{equation}\label{maximizationP}
\max_{Z\in\U\cap\Perp}\|\diag(Z^HAZ)\|_F^2.
\end{equation}

We develop the Jacobi-type algorithm for solving~\eqref{maximizationS} and~\eqref{maximizationP}.
In both cases this is an iterative algorithm
\begin{equation}\label{rel:jacobi}
A^{(k+1)}=R_k^HA^{(k)}R_k, \quad k\geq0, \quad A^{(0)}=A,
\end{equation}
where $R_k$ are structure-preserving rotations.
The goal of the $k$th step of~\eqref{rel:jacobi} is to make the Frobenius norm of the diagonal of $A^{(k+1)}$ as big as possible.
To achieve that we take the pivot pair $(i_k,j_k)$ and choose the appropriate rotation form and the rotation angles.
We obtain unitary structure-preserving matrix $Z$ that solves~\eqref{maximizationS} (or~\eqref{maximizationP}) as the product of these rotations. Then, we form the closest normal matrix as
$$X=Z(\diag(Z^HAZ)Z^H).$$

\subsection{Structure-preserving rotations}\label{sec:rotations}

Let us say more about the structure-preserving rotations used in~\eqref{rel:jacobi}
Symplectic and perplectic rotations that we use in each iterative step~\eqref{rel:jacobi} can be formed by embedding one or more Givens rotations
\begin{equation}\label{givens}
G=\left[
      \begin{array}{cc}
        \cos\phi & -e^{\imath\alpha}\sin\phi \\
        e^{-\imath\alpha}\sin\phi & \cos\phi \\
      \end{array}
    \right]
\end{equation}
into an identity matrix $I_{2n}$. To simplify the notation, set $c=\cos\phi$, $s=e^{\imath\alpha}\sin\phi$. Then matrix $G$ can be written as
$G=\left[
      \begin{array}{cc}
        c & -s \\
        \bar{s} & c \\
      \end{array}
    \right]$.
In our algorithm we use three kinds of symplectic and three kinds of perplectic embeddings to form rotations $R=R(i,j,\phi,\alpha)$.
Our rotations are similar to the structured rotations from~\cite{MMT03}, but note that the definitions of symplectic and perplectic matrices in~\cite{MMT03} slightly differ.

We start with symplectic rotations. If we insert only one Givens rotation $G$ from~\eqref{givens} into $I_{2n}$, we get a symplectic matrix only if this is done in a very special way. Matrix $G$ must be inserted on the intersection of $i$th and $(n+i)$th column and row and $\alpha$ must be zero. Therefore, we get
\begin{equation}\label{srot1}
R(i,j,\phi,\alpha)=R(i,n+i,\phi,0)=\left[
                         \begin{array}{ccccc}
                            &  &  &  &  \\
                            & \cos\phi &  & -\sin\phi &  \\
                            &  &  &  &  \\
                            & \sin\phi &  & \cos\phi &  \\
                            &  &  &  &  \\
                         \end{array}
                       \right]
\begin{array}{c}
   \\
  i \\
   \\
  n+i \\
   \\
\end{array}.
\end{equation}
All elements that are not explicitly written are as in $I_{2n}$.
Notice that in a Hamiltonian matrix entries on positions $(i,n+i)$ have only real and in skew-Hamiltonian matrix purely imaginary values. That is why real rotations are adequate here.
The second type of symplectic rotations that we use is the symplectic direct sum of two Givens rotations, that is,
\begin{equation}\label{srot2}
R(i,j,\phi,\alpha)={\footnotesize\left[
  \begin{array}{ccccc|ccccc}
     &  &  &  &  &  &  &  &  &  \\
     & c &  & -s &  &  &  &  &  &  \\
     &  &  &  &  &  &  &  &  &  \\
     & \bar{s} &  & c &  &  &  &  &  &  \\
     &  &  &  &  &  &  &  &  &  \\ \hline
     &  &  &  &  &  &  &  &  &  \\
     &  &  &  &  &  & c &  & -s &  \\
     &  &  &  &  &  &  &  &  &  \\
     &  &  &  &  &  & \bar{s} &  & c &  \\
     &  &  &  &  &  &  &  &  &  \\
  \end{array}
\right]
\begin{array}{c}
                         \\
                        i \\
                         \\
                        j \\
                         \\
                         \\
                        n+i \\
                         \\
                        n+j \\
                         \\
                      \end{array}}.
\end{equation}
Further on, we need concentric embedding of two Givens rotations given by
\begin{equation}\label{srot3}
R(i,j,\phi,\alpha)={\footnotesize\left[
  \begin{array}{ccccc|ccccc}
     &  &  &  &  &  &  &  &  &  \\
     & c &  &  &  &  &  &  & -s &  \\
     &  &  &  &  &  &  &  &  &  \\
     &  &  & c &  &  & -\bar{s} &  &  &  \\
     &  &  &  &  &  &  &  &  &  \\ \hline
     &  &  &  &  &  &  &  &  &  \\
     &  &  & s &  &  & c &  &  &  \\
     &  &  &  &  &  &  &  &  &  \\
     & \bar{s} &  &  &  &  &  &  & c &  \\
     &  &  &  &  &  &  &  &  &  \\
  \end{array}
\right]
\begin{array}{c}
                         \\
                        i \\
                         \\
                        j-n \\
                         \\
                         \\
                        n+i \\
                         \\
                        j \\
                         \\
                      \end{array}}.
\end{equation}

Pair $(i,j)$ in matrices~\eqref{srot1}, \eqref{srot2} and~\eqref{srot3} is called \emph{pivot pair}.
Usually in the Jacobi-type methods pivot pairs are taken from the upper triangle, $\mathcal{P}=\{(i,j) \ | \ 1\leq i<j\leq2n\}$.
Here, because of the double embeddings in rotations~\eqref{srot2} and~\eqref{srot3}, instead of $(i,j)$ one could equally say that the pivot pair is $(n+i,n+j)$, or $(j-n,n+i)$, respectively.
Therefore, we do not need to go through all pairs from $\mathcal{P}$, but its subset of $n^2$ positions.
\begin{align*}
\text{Rotations}~\eqref{srot1}: \quad & (i,n+i), \ 1\leq i\leq n & \leftarrow \ \text{$n$ pivot positions} & \\
\text{Rotations}~\eqref{srot2}: \quad & (i,j), \ 1\leq i<j\leq n & \leftarrow \ \text{$n(n-1)/2$ pivot positions} & \\
\text{Rotations}~\eqref{srot3}: \quad & (i,j), \ 1\leq i<n, n+i<j\leq n & \leftarrow \ \text{$n(n-1)/2$ pivot positions} &
\end{align*}

For better understanding we depict the pivot positions on a $10\times10$ matrix. Rotations~\eqref{srot1}, \eqref{srot2}, and~\eqref{srot3} are used on positions $\circ$, $\diamond$, and $\Box$, respectively,
\begin{equation}\label{pivotstrategy}
{\footnotesize
\left[
    \begin{array}{ccccc|ccccc}
       & \diamond & \diamond & \diamond & \diamond & \circ & \Box & \Box & \Box & \Box \\
       &  & \diamond & \diamond & \diamond &  & \circ & \Box & \Box & \Box \\
       &  &  & \diamond & \diamond &  &  & \circ & \Box & \Box \\
       &  &  &  & \diamond &  &  &  & \circ & \Box \\
       &  &  &  &  &  &  &  &  & \circ \\ \hline
       &  &  &  &  &  &  &  &  &  \\
       &  &  &  &  &  &  &  &  &  \\
       &  &  &  &  &  &  &  &  &  \\
       &  &  &  &  &  &  &  &  &  \\
       &  &  &  &  &  &  &  &  &  \\
    \end{array}
  \right]}.
\end{equation}
Considering the double embeddings on positions $\diamond$ and $\Box$, we see that the whole upper triangle is covered in the following way,
$${\footnotesize\left[
    \begin{array}{ccccc|ccccc}
       & \diamond & \diamond & \diamond & \diamond & \circ & \Box & \Box & \Box & \Box \\
       &  & \diamond & \diamond & \diamond & \Box & \circ & \Box & \Box & \Box \\
       &  &  & \diamond & \diamond & \Box & \Box & \circ & \Box & \Box \\
       &  &  &  & \diamond & \Box & \Box & \Box & \circ & \Box \\
       &  &  &  &  & \Box & \Box & \Box & \Box & \circ \\ \hline
       &  &  &  &  &  & \diamond & \diamond & \diamond & \diamond \\
       &  &  &  &  &  &  & \diamond & \diamond & \diamond \\
       &  &  &  &  &  &  &  & \diamond & \diamond \\
       &  &  &  &  &  &  &  &  & \diamond \\
       &  &  &  &  &  &  &  &  &  \\
    \end{array}
  \right]}.$$

On the other hand, perplectic rotation can be obtained by embedding only one Givens rotation $G$ only when $\alpha=-\frac{\pi}{2}$ and such $G$ is inserted on the intersection of the $i$th and the $(2n-i+1)$th column and row. Then we have
\begin{equation}\label{prot1}
R(i,j,\phi,\alpha)=R(i,2n-i+1,\phi,-\frac{\pi}{2})=\left[
                         \begin{array}{ccccc}
                            &  &  &  &  \\
                            & \cos\phi &  & \imath\sin\phi &  \\
                            &  &  &  &  \\
                            & \imath\sin\phi &  & \cos\phi &  \\
                            &  &  &  &  \\
                         \end{array}
                       \right]
\begin{array}{c}
   \\
  i \\
   \\
  2n-i+1 \\
   \\
\end{array}.
\end{equation}
When embedding two Givens rotations, we have more freedom. Perplectic direct sum embedding is given by
\begin{equation}\label{prot2}
R(i,j,\phi,\alpha)=
{\footnotesize\left[
  \begin{array}{ccccc|ccccc}
     &  &  &  &  &  &  &  &  &  \\
     & c &  & -s &  &  &  &  &  &  \\
     &  &  &  &  &  &  &  &  &  \\
     & \bar{s} &  & c &  &  &  &  &  &  \\
     &  &  &  &  &  &  &  &  &  \\ \hline
     &  &  &  &  &  &  &  &  &  \\
     &  &  &  &  &  & c &  & \bar{s} &  \\
     &  &  &  &  &  &  &  &  &  \\
     &  &  &  &  &  & -s &  & c &  \\
     &  &  &  &  &  &  &  &  &  \\
  \end{array}
\right]
\begin{array}{c}
                         \\
                        i \\
                         \\
                        j \\
                         \\
                         \\
                        2n-j+1 \\
                         \\
                        2n-i+1 \\
                         \\
                      \end{array}}.
\end{equation}
Finally, perplectic interleaved embedding of two Givens rotations is
\begin{equation}\label{prot3}
R(i,j,\phi,\alpha)=
{\footnotesize\left[
  \begin{array}{ccccc|ccccc}
     &  &  &  &  &  &  &  &  &  \\
     & c &  &  &  &  & -s &  &  &  \\
     &  &  &  &  &  &  &  &  &  \\
     &  &  & c &  &  &  &  & \bar{s} &  \\
     &  &  &  &  &  &  &  &  &  \\ \hline
     &  &  &  &  &  &  &  &  &  \\
     & \bar{s} &  &  &  &  & c &  &  &  \\
     &  &  &  &  &  &  &  &  &  \\
     &  &  & -s &  &  &  &  & c &  \\
     &  &  &  &  &  &  &  &  &  \\
  \end{array}
\right]
\begin{array}{c}
                         \\
                        i \\
                         \\
                        2n-j+1 \\
                         \\
                         \\
                        j \\
                         \\
                        2n-i+1 \\
                         \\
                      \end{array}}.
\end{equation}

Like it was the case with double embeddings~\eqref{srot2} and~\eqref{srot3}, for the pivot position in both rotations~\eqref{prot2} and~\eqref{prot3} one can also choose $(2n-j+1,2n-i+1)$ instead $(i,j)$. Here we are considering the following positions of pivot pairs.
\begin{align*}
\text{Rotations}~\eqref{prot1}: \quad & (i,j), \ 1\leq i\leq n, \ j=2n-i+1 & \leftarrow \ \text{$n$ pivot positions} & \\
\text{Rotations}~\eqref{prot2}: \quad & (i,j), \ 1\leq i<j\leq n & \leftarrow \ \text{$n(n-1)/2$ pivot positions} & \\
\text{Rotations}~\eqref{prot3}: \quad & (i,j), \ 1\leq i<n\leq j\leq2n-i & \leftarrow \ \text{$n(n-1)/2$ pivot positions} &
\end{align*}
Again, we depict this on a $10\times10$ matrix denoting the pivot positions corresponding to~\eqref{prot1}, \eqref{prot2}, and~\eqref{prot3} by $\circ$, $\diamond$, and $\Box$, respectively. We have
$${\footnotesize
\left[
    \begin{array}{ccccc|ccccc}
       & \diamond & \diamond & \diamond & \diamond & \Box & \Box & \Box & \Box & \circ \\
       &  & \diamond & \diamond & \diamond & \Box & \Box & \Box & \circ &  \\
       &  &  & \diamond & \diamond & \Box & \Box & \circ &  &  \\
       &  &  &  & \diamond & \Box & \circ &  &  &  \\
       &  &  &  &  & \circ &  &  &  &  \\ \hline
       &  &  &  &  &  &  &  &  &  \\
       &  &  &  &  &  &  &  &  &  \\
       &  &  &  &  &  &  &  &  &  \\
       &  &  &  &  &  &  &  &  &  \\
       &  &  &  &  &  &  &  &  &  \\
    \end{array}
  \right]}\overset{\text{considering double rotations}}{\rightsquigarrow}
{\footnotesize
\left[
    \begin{array}{ccccc|ccccc}
       & \diamond & \diamond & \diamond & \diamond & \Box & \Box & \Box & \Box & \circ \\
       &  & \diamond & \diamond & \diamond & \Box & \Box & \Box & \circ & \Box \\
       &  &  & \diamond & \diamond & \Box & \Box & \circ & \Box & \Box \\
       &  &  &  & \diamond & \Box & \circ & \Box & \Box & \Box \\
       &  &  &  &  & \circ & \Box & \Box & \Box & \Box \\ \hline
       &  &  &  &  &  & \diamond & \diamond & \diamond & \diamond \\
       &  &  &  &  &  &  & \diamond & \diamond & \diamond \\
       &  &  &  &  &  &  &  & \diamond & \diamond \\
       &  &  &  &  &  &  &  &  & \diamond \\
       &  &  &  &  &  &  &  &  &  \\
    \end{array}
  \right]}.$$

\subsection{Algorithm}

Knowing the shape of the structure-preserving rotations we write the Jacobi-type algorithm for solving the maximization problem~\eqref{maximizationS}. In each step we take a pivot pair $(i,j)$. The pivot position implies rotation form~\eqref{srot1}, \eqref{srot2} or~\eqref{srot3}, for which we compute the rotation angles. Then we perform one iterative step~\eqref{rel:jacobi} and update the symplectic unitary transformation matrix $Z_{k+1}=Z_kR_k$, $k\geq0$. Note that there is no need to set up matrices $R_k$.

\begin{Algorithm}\label{agm:Jacobi-Ham_closest}
\hrule\vspace{1ex}
\emph{Jacobi-type algorithm for solving maximization problem~\eqref{maximizationS}}
\vspace{0.5ex}\hrule
\begin{algorithmic}
\State \textbf{Input:} $A\in\mathbb{C}^{2n\times2n}$ Hamiltonian or skew-Hamiltonian.
\State \textbf{Output:} symplectic unitary $Z$
\State $k = 0$
\State $A^{(1)} = A$
\State $Z_1 = I$
\Repeat
\For{$i = 1, \ldots,n-1$}
\For{$j = i+1, \ldots, n$}
\State $k = k+1$
\State Find $\phi_k$, $\alpha_k$ for $R_k = R(i,j,\phi_k,\alpha_k)$ as in \eqref{srot2}
\State $A^{(k+1)}=R_k^HA^{(k)}R_k$
\State $Z_{k+1}=Z_kR_k$
\EndFor
\For{$j = n+i+1, \ldots, 2n$}
\State $k = k+1$
\State Find $\phi_k$, $\alpha_k$ for $R_k = R(i,j-n,\phi_k,\alpha_k)$ as in \eqref{srot3}
\State $A^{(k+1)}=R_k^HA^{(k)}R_k$
\State $Z_{k+1}=Z_kR_k$
\EndFor
\EndFor
\For{$i = 1, \ldots,n$}
\State $k = k+1$
\State Find $\phi_k$, $\alpha_k$ for $R_k = R(i,n+i,\phi_k,\alpha_k)$ as in \eqref{srot1}
\State $A^{(k+1)}=R_k^HA^{(k)}R_k$
\State $Z_{k+1}=Z_kR_k$
\EndFor
\Until{convergence}
\end{algorithmic}
\hrule
\end{Algorithm}

The order in which pivot pairs are taken defines a pivot strategy.
Algorithm~\ref{agm:Jacobi-Ham_closest} uses a cyclic pivot strategy. In general, cyclic strategies are periodic strategies with the period equal to the number of possible pivot positions. In our case that means that during the first $n^2$ steps in~\eqref{rel:jacobi} (and later on during any consecutive $n^2$ steps) we take all pivot positions corresponding to those marked in~\eqref{pivotstrategy}, each of them exactly once. This process defines one cycle. We repeat such cycles until the convergence is obtained. Specifically, reading from Algorithm~\ref{agm:Jacobi-Ham_closest}, we take pivot positions marked with $\diamond$ in the row-wise order, then positions $\Box$ in the row-wise order, and positions $\circ$ again row by row.

As we will see in Section~\ref{sec:cvgproof}, the convergence does not depend on the order inside one cycle and our proof holds for any cyclic pivot strategy. Therefore, \textit{for} loops in Algorithm~\ref{agm:Jacobi-Ham_closest} can be altered depending on the pivot strategy. Nevertheless, in order to ensure the convergence one must check that each pivot pair satisfies the condition of Lemma~\ref{lemma:pivotcondition-Ham},
$$|\langle \emph{grad} f_{\Ham}(Z),Z\dot{R}(i,j,0,\alpha) \rangle|\geq\frac{2}{\sqrt{4n^2-2n}}\|\emph{grad} f_{\Ham}(Z)\|_F.$$
If this is not true for some pivot pair, that pair is skipped.

Algorithm~\ref{agm:Jacobi-Ham_closest} can easily be modified for solving minimization problem~\eqref{maximizationP}.
Instead of rotations~\eqref{srot1}, \eqref{srot2} and~\eqref{srot3} perplectic rotations~\eqref{prot1}, \eqref{prot2} and~\eqref{prot3} are used.
The inner \textit{for} loop for the rotation~\eqref{prot2} has to be modified to 'for $j=n+1:2n-i$ do'.

\subsection{The choice of rotation angles $\phi$ and $\alpha$}

In the $k$th step of Algorithm~\ref{agm:Jacobi-Ham_closest} one should chose rotation angles $\phi_k$ and $\alpha_k$ such that
$R_k=R(i_k,j_k,\phi_k,\alpha_k)$ maximizes the Frobenius norm of the diagonal of $A^{(k+1)}=R_k^HA^{(k)}R_k$.
Here we show how $\phi$ and $\alpha$ are obtained in the case of symplectic rotations. The same reasoning holds for perplectic rotations.

Denote the $k$th pivot position by $(i_k,i_k)$. If $R_k$ is of the form~\eqref{srot1}, after one iteration two diagonal elements on positions $(i_k,i_k)$ and $(j_k,j_k)$ are changed. If $R_k$ is double Givens rotation~\eqref{srot2} or~\eqref{srot3}, then four diagonal elements are changed. However, for $A$ Hamiltonian or skew-Hamiltonian matrix we have
\begin{subequations}
\begin{align}
|\re(a_{n+i,n+i})| & =|\re(a_{ii})|, & |\im(a_{n+i,n+i})| & =|\im(a_{ii})|, \quad 1\leq i\leq n \label{angle:eq1}, \\
|\re(a_{n+j,n+j})| & =|\re(a_{jj})|, & |\im(a_{n+j,n+j})| & =|\im(a_{jj})|, \quad 1\leq j\leq n \label{angle:eq2}, \\
|\re(a_{j-n,j-n})| & =|\re(a_{jj})|, & |\im(a_{j-n,j-n})| & =|\im(a_{jj})|, \quad n+i\leq j\leq2n \label{angle:eq3}.
\end{align}
\end{subequations}
From~\eqref{angle:eq1} and~\eqref{angle:eq2} it follows that for rotations~\eqref{srot2} it is enough to consider only the changes on positions $(i_k,i_k)$ and $(j_k,j_k)$, $1\leq i_k,j_k\leq n$. The same conclusion follows from~\eqref{angle:eq1} and~\eqref{angle:eq3} for rotations~\eqref{srot3}, with $1\leq i_k\leq n<j_k\leq2n$.
Thus it is always enough to consider only the changes induced by one Givens rotation and the following computation holds for all rotations~\eqref{srot1}, \eqref{srot2}, and~\eqref{srot3}.

For the simplicity of notation, denote $A^{(k+1)}=A'=(a_{ij}')$, $A^{(k)}=A=(a_{ij})$, $\phi_k=\phi$, $\alpha_k=\alpha$.
Consider the pivot submatrix
\begin{equation}\label{angle:submatrix}
\left[
    \begin{array}{cc}
      a_{ii}' & a_{ij}' \\
      a_{ji}' & a_{jj}' \\
    \end{array}
  \right] = \left[
    \begin{array}{cc}
      \cos\phi & -e^{\imath\alpha}\sin\phi \\
      e^{-\imath\alpha}\sin\phi & \cos\phi \\
    \end{array}
  \right]^H\left[
    \begin{array}{cc}
      a_{ii} & a_{ij} \\
      a_{ji} & a_{jj} \\
    \end{array}
  \right]\left[
    \begin{array}{cc}
      \cos\phi & -e^{\imath\alpha}\sin\phi \\
      e^{-\imath\alpha}\sin\phi & \cos\phi \\
    \end{array}
  \right].
\end{equation}
We need
\begin{equation}\label{angle:max1}
|a_{ii}'|^2+|a_{jj}'|^2\rightarrow\max.
\end{equation}

Set $a_{rs}=x_{rs}+y_{rs}\imath$. Condition~\eqref{angle:max1} becomes
$$|x_{ii}'+y_{ii}'\imath|^2+|x_{jj}'+y_{jj}'\imath|^2 =
(x_{ii}')^2+(y_{ii}')^2+(x_{jj}')^2+(y_{jj}')^2 \rightarrow\max.$$
From~\eqref{angle:submatrix} it follows
\begin{align*}
x_{ii}'+y_{ii}'\imath & =(x_{ii}+y_{ii}\imath)\cos^2\phi+(x_{jj}+y_{jj}\imath)\sin^2\phi \\
& \quad +(x_{ij}\cos\alpha-x_{ij}\sin\alpha\imath+y_{ij}\cos\alpha\imath+y_{ij}\sin\alpha)\sin\phi\cos\phi \\
& \quad +(x_{ji}\cos\alpha+x_{ji}\sin\alpha\imath+y_{ji}\cos\alpha\imath-y_{ji}\sin\alpha)\sin\phi\cos\phi, \\
x_{jj}'+y_{jj}'\imath & =(x_{ii}+y_{ii}\imath)\sin^2\phi+(x_{jj}+y_{jj}\imath)\cos^2\phi \\
& \quad -(x_{ij}\cos\alpha-x_{ij}\sin\alpha\imath+y_{ij}\cos\alpha\imath+y_{ij}\sin\alpha)\sin\phi\cos\phi \\
& \quad -(x_{ji}\cos\alpha+x_{ji}\sin\alpha\imath+y_{ji}\cos\alpha\imath-y_{ji}\sin\alpha)\sin\phi\cos\phi.
\end{align*}
Splitting the real and imaginary part and using $e^{\imath\alpha}=\cos\alpha+\imath\sin\alpha$ gives
\begin{align*}
x_{ii}' & =x_{ii}\cos^2\phi+x_{jj}\sin^2\phi
+(x_{ij}\cos\alpha+y_{ij}\sin\alpha)\sin\phi\cos\phi +(x_{ji}\cos\alpha-y_{ji}\sin\alpha)\sin\phi\cos\phi, \\
y_{ii}' & =y_{ii}\cos^2\phi+y_{jj}\sin^2\phi
+(-x_{ij}\sin\alpha+y_{ij}\cos\alpha)\sin\phi\cos\phi +(x_{ji}\sin\alpha+y_{ji}\cos\alpha)\sin\phi\cos\phi, \\
x_{jj}' & =x_{ii}\sin^2\phi+x_{jj}\cos^2\phi
-(x_{ij}\cos\alpha+y_{ij}\sin\alpha)\sin\phi\cos\phi -(x_{ji}\cos\alpha-y_{ji}\sin\alpha)\sin\phi\cos\phi, \\
y_{jj}' & =y_{ii}\sin^2\phi+y_{jj}\cos^2\phi
-(-x_{ij}\sin\alpha+y_{ij}\cos\alpha)\sin\phi\cos\phi -(x_{ji}\sin\alpha+y_{ji}\cos\alpha)\sin\phi\cos\phi.
\end{align*}

Define the function $g:\langle-\frac{\pi}{4},\frac{\pi}{4}]\times\langle-\frac{\pi}{2},\frac{\pi}{2}]\rightarrow\mathbb{R}$,
\begin{equation}\label{angle:function}
g(\phi,\alpha)=(x_{ii}')^2+(y_{ii}')^2+(x_{jj}')^2+(y_{jj}')^2.
\end{equation}
Finding $\phi$ and $\alpha$ that maximize $g$ will give the solution of the maximization problem~\eqref{angle:max1}.
Partial derivatives of $g$ are
\begin{align}
0=\frac{\partial}{\partial\phi}g(\phi,\alpha) & =
2\cos\alpha\cos4\phi\big((x_{ij}+x_{ji})(x_{ii}-x_{jj})+(y_{ij}+y_{ji})(y_{ii}-y_{jj})\big) \nonumber \\
& \quad +2\sin\alpha\cos4\phi\big((x_{ii}-x_{jj})(y_{ij}-y_{ji})+(x_{ji}-x_{ij})(y_{ii}-y_{jj})\big) \nonumber \\
& \quad +\sin4\phi\big(x_{ij}^2+x_{ji}^2+y_{ij}^2+y_{ji}^2-(x_{ii}-x_{jj})^2-(y_{ii}-y_{jj})^2 \nonumber \\
& \quad +2\cos2\alpha(x_{ij}x_{ji}+y_{ij}y_{ji}) + 2\sin2\alpha(x_{ji}y_{ij}-x_{ij}y_{ji})\big), \label{angle:dphi} \\
0=\frac{\partial}{\partial\alpha}g(\phi,\alpha) & =
2\sin^22\phi\big((x_{ji}y_{ij}-x_{ij}y_{ji})\cos2\alpha-(x_{ij}x_{ji}+y_{ij}y_{ji})\sin2\alpha\big) \nonumber \\
& \quad +\sin2\phi\cos2\phi\Big(\big((x_{ii}-x_{jj})(y_{ij}-y_{ji})-(x_{ij}-x_{ji})(y_{ii}-y_{jj})\big)\cos\alpha \nonumber \\
& \quad +\big((x_{ij}+x_{ji})(x_{jj}-x_{ii})-(y_{ij}+y_{ji})(y_{ii}-y_{jj})\big)\sin\alpha\Big). \label{angle:dalpha}
\end{align}

We take a closer look at~\eqref{angle:dalpha}
and distinguish between different cases.
\begin{itemize}
\item The trivial solution is $\phi=0$. The transformation matrix $R$ will be the identity.
\item If $\phi=\frac{\pi}{4}$, relation~\eqref{angle:dalpha} simplifies to
$$0=(x_{ji}y_{ij}-x_{ij}y_{ji})\cos2\alpha-(x_{ij}x_{ji}+y_{ij}y_{ji})\sin2\alpha.$$
Then we either have $x_{ij}x_{ji}+y_{ij}y_{ji}=0$ and $\alpha=\pm\frac{\pi}{4}$ or
\begin{equation}\label{angle:t2alpha}
\tan2\alpha=\frac{x_{ji}y_{ij}-x_{ij}y_{ji}}{x_{ij}x_{ji}+y_{ij}y_{ji}}.
\end{equation}
\end{itemize}

Otherwise, we divide~\eqref{angle:dalpha} by $\cos^22\phi$ and set $t=\tan2\phi$. We obtain a quadratic equation in $t$,
$$K_2(\alpha)t^2+K_1(\alpha)t=0,$$
where
\begin{align*}
K_1(\alpha) & = \big((x_{ii}-x_{jj})(y_{ij}-y_{ji})-(x_{ij}-x_{ji})(y_{ii}-y_{jj})\big)\cos\alpha \\
& \quad +\big((x_{ij}+x_{ji})(x_{jj}-x_{ii})-(y_{ij}+y_{ji})(y_{ii}-y_{jj})\big)\sin\alpha, \\
K_2(\alpha) & = 2\big((x_{ji}y_{ij}-x_{ij}y_{ji})\cos2\alpha-(x_{ij}x_{ji}+y_{ij}y_{ji})\sin2\alpha\big).
\end{align*}
Since $\phi\neq0$, we have $t\neq0$ and
\begin{equation}\label{angle:phi}
t=-\frac{K_1(\alpha)}{K_2(\alpha)}.
\end{equation}
We have $K_2(\alpha)\neq0$ because $\phi\neq\frac{\pi}{4}$.

Now we consider relation~\eqref{angle:dphi}.
Again, we distinguish between different cases.
\begin{itemize}
\item If $\alpha=\frac{\pi}{2}$, relation~\eqref{angle:dphi} simplifies to
\begin{align*}
0 & = 2\big((x_{ii}-x_{jj})(y_{ij}-y_{ji})+(x_{ji}-x_{ij})(y_{ii}-y_{jj})\big)\cos4\phi \\
& \quad +\big((x_{ij}-x_{ji})^2+(y_{ij}-y_{ji})^2-(x_{ii}-x_{jj})^2-(y_{ii}-y_{jj})^2\big)\sin4\phi.
\end{align*}
Then we either have $\phi=\pm\frac{\pi}{8}$ or
\begin{equation}\label{angle:t4phi}
\tan4\phi=\frac{-2\big((x_{ii}-x_{jj})(y_{ij}-y_{ji})+(x_{ji}-x_{ij})(y_{ii}-y_{jj})\big)}
{(x_{ij}-x_{ji})^2+(y_{ij}-y_{ji})^2-(x_{ii}-x_{jj})^2-(y_{ii}-y_{jj})^2}.
\end{equation}
\end{itemize}

Otherwise, we substitute
\begin{align*}
\cos4\phi & = \frac{1-t^2}{1+t^2}, & \sin4\phi & = \frac{2t}{1+t^2}, \\
\cos2\alpha & = \cos^2\alpha-\sin^2\alpha, & \sin2\alpha & = 2\sin\alpha\cos\alpha
\end{align*}
in~\eqref{angle:dphi} and obtain
\begin{align*}
&2\cos\alpha(1-t^2)\big((x_{ij}+x_{ji})(x_{ii}-x_{jj})+(y_{ij}+y_{ji})(y_{ii}-y_{jj})\big) \nonumber \\
& \quad +2\sin\alpha(1-t^2)\big((x_{ii}-x_{jj})(y_{ij}-y_{ji})+(x_{ji}-x_{ij})(y_{ii}-y_{jj})\big) \nonumber \\
& \quad +2t\big(x_{ij}^2+x_{ji}^2+y_{ij}^2+y_{ji}^2-(x_{ii}-x_{jj})^2-(y_{ii}-y_{jj})^2 \nonumber \\
& \quad +2(\cos^2\alpha-\sin^2\alpha)(x_{ij}x_{ji}+y_{ij}y_{ji}) + 4\cos\alpha\sin\alpha(x_{ji}y_{ij}-x_{ij}y_{ji})\big)=0.
\end{align*}
Using~\eqref{angle:phi} we multiply the obtained equation with $K_2(\alpha)^2$ and get
\begin{align}
&\cos\alpha(K_2(\alpha)^2-K_1(\alpha)^2)\big((x_{ij}+x_{ji})(x_{ii}-x_{jj})+(y_{ij}+y_{ji})(y_{ii}-y_{jj})\big) \nonumber \\
& \quad +\sin\alpha(K_2(\alpha)^2-K_1(\alpha)^2)\big((x_{ii}-x_{jj})(y_{ij}-y_{ji})+(x_{ji}-x_{ij})(y_{ii}-y_{jj})\big) \nonumber \\
& \quad -K_1(\alpha)K_2(\alpha)\big(x_{ij}^2+x_{ji}^2+y_{ij}^2+y_{ji}^2-(x_{ii}-x_{jj})^2-(y_{ii}-y_{jj})^2 \nonumber \\
& \quad +2(\cos^2\alpha-\sin^2\alpha)(x_{ij}x_{ji}+y_{ij}y_{ji}) + 4\cos\alpha\sin\alpha(x_{ji}y_{ij}-x_{ij}y_{ji})\big)=0. \label{angle:eq5}
\end{align}

The left-hand side in~\eqref{angle:eq5} is a sum of expressions of the form $C\cos^k\alpha\sin^l\alpha$, for $k+l=3$ or $k+l=5$, and different $C\in\mathbb{R}$.
If we take a closer look, we see that the summands where $k+l=5$ can be reduced. Precisely, the sum of all expressions $C\cos^k\alpha\sin^l\alpha$ such that $k+l=5$ equals
\begin{align*}
& -8 \Big((x_{ij}^2 x_{ji} y_{ii} - x_{ij} x_{ji}^2 y_{ii} + x_{ii} x_{ji}^2 y_{ij} -
x_{ji}^2 x_{jj} y_{ij}  + x_{ji} y_{ii} y_{ij}^2  - x_{ii} x_{ij}^2 y_{ji}  +
x_{ij}^2 x_{jj} y_{ji}  - x_{ii} y_{ij}^2 y_{ji} \\
& \quad + x_{jj} y_{ij}^2 y_{ji}  -
x_{ij} y_{ii} y_{ji}^2  + x_{ii} y_{ij} y_{ji}^2  - x_{jj} y_{ij} y_{ji}^2  -
x_{ij}^2 x_{ji} y_{jj}  + x_{ij} x_{ji}^2 y_{jj}  - x_{ji} y_{ij}^2 y_{jj}  +
x_{ij} y_{ji}^2 y_{jj})\cos\alpha \\
& \quad +(-x_{ii} x_{ij}^2 x_{ji} - x_{ii} x_{ij} x_{ji}^2  +
x_{ij}^2 x_{ji} x_{jj}  + x_{ij} x_{ji}^2 x_{jj}  - x_{ji}^2 y_{ii} y_{ij}  -
x_{ii} x_{ji} y_{ij}^2  + x_{ji} x_{jj} y_{ij}^2  - x_{ij}^2 y_{ii} y_{ji} \\
& \quad -y_{ii} y_{ij}^2 y_{ji}  - x_{ii} x_{ij} y_{ji}^2  + x_{ij} x_{jj} y_{ji}^2  -
y_{ii} y_{ij} y_{ji}^2  + x_{ji}^2 y_{ij} y_{jj}  + x_{ij}^2 y_{ji} y_{jj}  +
y_{ij}^2 y_{ji} y_{jj}  + y_{ij} y_{ji}^2 y_{jj}) \sin\alpha\Big) \cdot \\
& \cdot (\cos^2\alpha + \sin^2\alpha) \cdot \\
& \cdot (-x_{ji} y_{ij} \cos^2\alpha + x_{ij} y_{ji} \cos^2\alpha + 2x_{ij} x_{ji} \cos\alpha \sin\alpha +
2 y_{ij} y_{ji} \cos\alpha \sin\alpha + x_{ji} y_{ij} \sin^2\alpha - x_{ij} y_{ji} \sin^2\alpha).
\end{align*}
Using only the fact that $\cos^2\alpha + \sin^2\alpha=1$ we see that this is again a sum of expressions $C\cos^k\alpha\sin^l\alpha$ for $k+l=3$.
Therefore, the left-hand side in~\eqref{angle:eq5} is a sum of expressions of the form $C\cos^k\alpha\sin^l\alpha$ for $k+l=3$.

Thus we can divide equation~\eqref{angle:eq5} by $\cos^3\alpha$ and set $\tau=\tan\alpha$.
We get a cubic equation in $\tau$,
\begin{equation}\label{angle:C}
C_3\tau^3+C_2\tau^2+C_1+C_0\tau=0.
\end{equation}
Equation~\eqref{angle:C} has at least one real solution. For each real solution we substitute
$\alpha=\arctan\tau$ into~\eqref{angle:phi} to obtain $t$, and hence $\phi =\frac{1}{2}\arctan t$.

Finally, we take the pair $(\phi,\alpha)$ from among all possible solutions that gives the largest value of function $g$.
Algorithm~\ref{agm:angle} summarizes the process of computing $\phi$ and $\alpha$.

\begin{Algorithm}\label{agm:angle}
\hrule\vspace{1ex}
\emph{Rotation angles in Algorithm \ref{agm:Jacobi-Ham_closest}}
\vspace{0.5ex}\hrule
\begin{algorithmic}
\State Form $g(\phi,\alpha)$ as in \eqref{angle:function}.
\State Case 1: $(\phi_1,\alpha_1) = (0, 0)$.
\State Case 2: $(\phi_2,\alpha_2) = (\frac{\pi}{4},\pm\frac{\pi}{4})$, or
$(\phi_2,\alpha_2) = (\frac{\pi}{4},\alpha_2)$
with $\alpha_2$ as in~\eqref{angle:t2alpha}.
\State Case 3: $(\phi_3,\alpha_3) = (\pm\frac{\pi}{8},\frac{\pi}{2})$, or
$(\phi_3,\alpha_3) = (\phi_3,\frac{\pi}{2})$
with $\phi$ as in~\eqref{angle:t4phi}.
\State Case 4: $(\phi_4,\alpha_4)$ with $\alpha_4=\arctan\tau$ for all real solutions $\tau$ of~\eqref{angle:C} and the corresponding $\phi_4$ from~\eqref{angle:phi}.
\State Choose that pair $(\phi,\alpha)$ which gives the largest value of $g(\phi,\alpha)$.
\end{algorithmic}
\hrule
\end{Algorithm}

\section{Convergence of Algorithm~\ref{agm:Jacobi-Ham_closest}}\label{sec:cvgproof}

In this section we provide a convergence proof for Algorithm~\ref{agm:Jacobi-Ham_closest}. We will discuss only the case of
Hamiltonian matrices, the proof for the other three structures follows in the same way with only minor modifications.
Related to the maximization problem~\eqref{maximizationS} we define the objective function
\begin{equation}\label{function}
f_\Ham : \U\cap\Symp\rightarrow\mathbb{R}_{\geq0}, \quad
f_\Ham(Z) =\|\diag(Z^HHZ)\|_F^2=\sum_{j=1}^{2n}|\langle AZe_j,Ze_j \rangle|^2.
\end{equation}
We show that Algorithm~\ref{agm:Jacobi-Ham_closest} converges to the stationary point of this function.
In particular, we will prove the following theorem.

\begin{Theorem}\label{tm:cvg-Ham}
Let $(Z_k,\ k\geq0)$ be the sequence generated by Algorithm~\ref{agm:Jacobi-Ham_closest}. Every accumulation point of $(Z_k,\ k\geq0)$ is a stationary point of function $f_{\Ham}$ from~\eqref{function}.
\end{Theorem}

The proof of Theorem~\ref{tm:cvg-Ham} uses the technique from~\cite{Ishteva13}.
and is based on Polak's theorem on model algorithms~\cite[Section 1.3, Theorem 3]{Polak71}. We will need three auxiliary results from Lemmas~\ref{lemma:gradf-Ham}, \ref{lemma:pivotcondition-Ham}, and~\eqref{lemma:cvgdelta-Ham}.

Before we move to Lemma~\ref{lemma:gradf-Ham} that gives the structure of $\emph{grad}f_{\Ham}(Z)$, let us say a bit more about the function $f_{\Ham}$.
As a real valued function of a complex variable is complex differentiable only if it is constant, the function $f_\Ham$ is not complex
differentiable. That is, $\frac{\partial f_\Ham}{\partial z_{jk}}$ does not exist for $Z=[z_{jk}]\in\mathbb{C}^{2n\times2n}$.
But with $z_{jk} = \re(z_{jk})+\imath\im(z_{jk})$, partial derivatives
$$\frac{\partial f_\Ham}{\partial \re(z_{jk})} \textnormal{\quad and \quad} \frac{\partial f_\Ham}{\partial\im(z_{jk})}, \quad j,k=1,\ldots,2n,$$
do exist as this involves only real differentiation. We identify $\mathbb{C}$ with $\mathbb{R}^{1\times2}$,
$\mathbb{C}^{2n}$ with $\mathbb{R}^{2n\times 2}$ and $\mathbb{C}^{2n \times 2n}$ with $\mathbb{R}^{2n\times2n\times2}$.
Consequently, $f_\Ham$ is viewed as a real-valued function on the Euclidian space $\mathbb{R}^{2n\times2n\times2}$.
As such, it is differentiable and its gradient is a matrix
$$\grad f_\Ham(Z) = \left[\frac{\partial f_\Ham}{\partial\re(z_{jk})} + \imath \frac{\partial f_\Ham}{\partial\im(z_{jk})}\right]_{j,k=1}^{2n}.$$

\begin{Lemma}\label{lemma:gradf-Ham}
The gradient of $f_{\Ham}$ from~\eqref{function} can be expressed as
$$\emph{grad}f_{\Ham}(Z)=ZX,$$
where $\diag(X)=0$ and $X$ is skew-Hermitian Hamiltonian.
\end{Lemma}

\begin{proof}
We will not be able to determine $\grad f_\Ham(Z)$ directly. Instead, we define function
$\widetilde{f}:\mathbb{C}^{2n\times2n}\rightarrow\mathbb{R}_{\geq 0}$,
$$\widetilde{f}(Z)=\sum_{j=1}^{2n}|\langle AZe_j,Ze_j \rangle|^2,$$
on a larger domain. Then $f_\Ham$ is the restriction of $\widetilde{f}$ to $\U\cap\Symp$. We first determine
$\grad \widetilde{f}(Z).$

To that end we define a new function $g:\mathbb{C}^{2n}\rightarrow\mathbb{R}$,
$g(z)=|\langle Az,z \rangle|^2$.
It allows us to rewrite $\widetilde{f}$ as
$$\widetilde{f}(Z)=\sum_{j=1}^{2n} g(Ze_j).$$
Then
$$\grad\tilde{f}(Z)=\left[
                  \begin{array}{ccc}
                    \nabla g(Ze_1) & \cdots & \nabla g(Ze_{2n}) \\
                  \end{array}
                \right],$$
where
$$\nabla g(z) = \left[ \frac{\partial g}{\partial \re(z_j)} + \imath \frac{\partial g}{\partial \im(z_j)} \right]_{j=1}^{2n},$$
as $g$ is real differentiable.

In order to determine $\nabla g(z)$ we use Taylor expansion of $g$,
\begin{equation}\label{rel:taylor}
g(z+h)=g(z)+\langle\nabla g(z),h\rangle_\mathbb{R} + O(\|h\|^2),
\end{equation}
for $h \in \mathbb{C}^{2n}$ and
$\langle u,v \rangle_\mathbb{R}= \re(\langle u,v \rangle)$.
We have
\begin{align*}
g(z+h)-g(z) & = |\langle A(z+h),z+h\rangle|^2-|\langle Az,z\rangle|^2 \\
& = |\langle Az,z\rangle + \langle Az,h\rangle + \langle Ah,z\rangle + \langle Ah,h\rangle|^2-|\langle Az,z\rangle|^2 \\
& = 2\re((\langle Az,h\rangle + \langle Ah,z\rangle)\overline{\langle Az,z\rangle}) + O(\|h\|^2) \\
& = \re(\langle 2\overline{\langle Az,z\rangle} Az + 2\langle Az,z\rangle A^Hz, h\rangle) + O(\|h\|^2).
\end{align*}
Then
$$
g(z+h)-g(z)=2\big\langle\overline{\langle Az,z\rangle} Az + \langle Az,z\rangle A^Hz, h\big\rangle_\mathbb{R} + O(\|h\|^2).$$
Relation~\eqref{rel:taylor} implies
$$\nabla g(z)=2\overline{\langle Az,z\rangle} Az + 2\langle Az,z\rangle A^Hz.$$
With this, we have described $\grad\widetilde{f}(Z).$

Further on, $\grad f_\Ham(Z)$ is obtained by projecting $\grad\widetilde{f}(Z)$ onto the tangent space of unitary symplectic matrices at $Z$.
We have
$$\grad f_\Ham(Z) = \pi(\grad \widetilde{f}(Z)).$$
For any unitary (and symplectic) matrix $Z$, matrix
\begin{equation}\label{eq:lemmaY}
Y:=Z^H\grad\widetilde{f}(Z)
\end{equation}
does exist. Thus we can write
$\grad\widetilde{f}(Z)=ZY$.
Then
$$\grad f_\Ham(Z) = ZX,$$
where $X=\pi(Y)$.
Since the tangent space of the group of unitary symplectic matrices at the identity are skew-Hermitian Hamiltonian matrices, it only remains to prove that $\diag(X)=0$.

The diagonal of $Y$ from~\eqref{eq:lemmaY} is given by
$$\diag(Y)=\diag(Z^H \grad\tilde{f}(Z))=(\langle\nabla g(Ze_j),Ze_j\rangle)_{j=1}^{2n}.$$
Further on,
$$\langle\nabla g(z),z\rangle = \langle 2\overline{\langle Az,z\rangle} Az + 2\langle Az,z\rangle A^Hz = 4|\langle Az,z\rangle|^2 \in\mathbb{R}.$$
Therefore, $\diag(Y)$ is real. Its projection onto the space of skew-Hermitian matrices will give zeros on the diagonal of $X$, that is $\diag(X)=0$.
\end{proof}

\begin{Remark}
The orthogonal projection of $Y=\left[
     \begin{smallmatrix}
       Y_{11} & Y_{12} \\
       Y_{21} & Y_{22}
     \end{smallmatrix}
   \right]
$ onto the subspace of skew-Hermitian Hamiltonian matrices is given by
$$\left[
     \begin{array}{cc}
       B & C \\
       -C & B \\
     \end{array}
   \right], \qquad B=\frac{Y_{11}+Y_{22}-Y_{11}^H-Y_{22}^H}{4}, \quad C=\frac{Y_{12}-Y_{21}+Y_{12}^H-Y_{21}^H}{4}.$$
\end{Remark}

Lemma~\ref{lemma:gradf-Ham} is used in Lemma~\ref{lemma:pivotcondition-Ham}.

\begin{Lemma}\label{lemma:pivotcondition-Ham}
For every symplectic unitary $Z\in\mathbb{C}^{2n\times2n}$ there is symplectic rotation $R(i,j,\phi,\alpha)$ such that
$$|\langle \emph{grad} f_{\Ham}(Z),Z\dot{R}(i,j,0,\alpha) \rangle|\geq\eta\|\emph{grad} f_{\Ham}(Z)\|_F, \qquad \eta=\frac{2}{\sqrt{4n^2-2n}},$$
where $\dot{R}(i,j,0,\alpha)$ denotes $\frac{\partial}{\partial\phi}R(i,j,\phi,\alpha)\Big{|}_{\phi=0}$.
\end{Lemma}

\begin{proof}
Obviously, if $\|\grad f_{\Ham}(Z)\|_F=0$, the assertion holds for any rotation. Thus, assume that $\|\grad f_{\Ham}(Z)\|_F\neq0$.
From Lemma~\ref{lemma:gradf-Ham} we know that $\grad f_{\Ham}(Z)=ZX.$ Hence, for $X=[x_{ij}]_{i,j=1}^{2n}$ and
$|x|=|x_{rs}|=\max_{i\neq j}|x_{ij}|>0$ we have
\begin{equation}\label{pivot-gradfmax}
\|\grad f_{\Ham}(Z)\|_F=\|ZX\|_F=\|X\|_F\leq\sqrt{4n^2-2n} \ |x|.
\end{equation}
On the other hand,
\begin{align}
\langle \grad f_{\Ham}(Z),Z\dot{R}(i,j,0,\alpha) \rangle_\mathbb{R} &
= \re(\tr((\grad f_{\Ham}(Z))^HZ\dot{R}(i,j,0,\alpha))) \nonumber \\
& = \re(\tr((ZX)^HZ\dot{R}(i,j,0,\alpha)))  \nonumber \\
& = \re(\tr(X^H\dot{R}(i,j,0,\alpha))). \label{pivot-trace}
\end{align}

Let us first consider a unitary symplectic rotation $R(i,j,\phi,\alpha)$ of the form~\eqref{srot2}. Then
$$\dot{R}(i,j,\phi,\alpha)={\footnotesize\left[
  \begin{array}{ccccc|ccccc}
     &  &  &  &  &  &  &  &  &  \\
     & -\sin\phi &  & -e^{\imath\alpha}\cos\phi &  &  &  &  &  &  \\
     &  &  &  &  &  &  &  &  &  \\
     & e^{-\imath\alpha}\cos\phi &  & -\sin\phi &  &  &  &  &  &  \\
     &  &  &  &  &  &  &  &  &  \\ \hline
     &  &  &  &  &  &  &  &  &  \\
     &  &  &  &  &  & -\sin\phi &  & -e^{\imath\alpha}\cos\phi &  \\
     &  &  &  &  &  &  &  &  &  \\
     &  &  &  &  &  & e^{-\imath\alpha}\cos\phi &  & -\sin\phi &  \\
     &  &  &  &  &  &  &  &  &  \\
  \end{array}
\right]}$$
and
$$\dot{R}(i,j,0,\alpha)={\footnotesize\left[
  \begin{array}{ccccc|ccccc}
     &  &  &  &  &  &  &  &  &  \\
     &  &  & -e^{\imath\alpha} &  &  &  &  &  &  \\
     &  &  &  &  &  &  &  &  &  \\
     & e^{-\imath\alpha} &  &  &  &  &  &  &  &  \\
     &  &  &  &  &  &  &  &  &  \\ \hline
     &  &  &  &  &  &  &  &  &  \\
     &  &  &  &  &  &  &  & -e^{\imath\alpha} &  \\
     &  &  &  &  &  &  &  &  &  \\
     &  &  &  &  &  & e^{-\imath\alpha} &  &  &  \\
     &  &  &  &  &  &  &  &  &  \\
  \end{array}
\right]},$$
where in the matrices on the right-hand side all elements that are not explicitly given are zero.
It follows from Lemma~\ref{lemma:gradf-Ham} that matrix $X$ is skew-Hermitian and Hamiltonian. For $x = x_{ij}$ we have $x_{ji}=-\bar{x}$, $x_{n+i,n+j}=x$ and $x_{n+j,n+i}=-\bar{x}$. This gives
$$X^H\dot{R}(i,j,0,\alpha)
= {\footnotesize\left[
  \begin{array}{ccccc|ccccc}
     &  &  &  &  &  &  &  &  &  \\
     &  &  & -x &  &  &  &  &  &  \\
     &  &  &  &  &  &  &  &  &  \\
     & \bar{x} &  &  &  &  &  &  &  &  \\
     &  &  &  &  &  &  &  &  &  \\ \hline
     &  &  &  &  &  &  &  &  &  \\
     &  &  &  &  &  &  &  & -x &  \\
     &  &  &  &  &  &  &  &  &  \\
     &  &  &  &  &  & \bar{x} &  &  &  \\
     &  &  &  &  &  &  &  &  &  \\
  \end{array}
\right]\left[
  \begin{array}{ccccc|ccccc}
     &  &  &  &  &  &  &  &  &  \\
     &  &  & -e^{\imath\alpha} &  &  &  &  &  &  \\
     &  &  &  &  &  &  &  &  &  \\
     & e^{-\imath\alpha} &  &  &  &  &  &  &  &  \\
     &  &  &  &  &  &  &  &  &  \\ \hline
     &  &  &  &  &  &  &  &  &  \\
     &  &  &  &  &  &  &  & -e^{\imath\alpha} &  \\
     &  &  &  &  &  &  &  &  &  \\
     &  &  &  &  &  & e^{-\imath\alpha} &  &  &  \\
     &  &  &  &  &  &  &  &  &  \\
  \end{array}
\right]},$$
where in $X^H$ only four relevant entries at the positions $(i,j)$, $(j,i)$, $(n+i,n+j)$ and $(n+j,n+i)$ are given.
Now~\eqref{pivot-trace} implies
$$\langle \grad f_{\Ham}(Z),Z\dot{R}(i,j,0,\alpha) \rangle_\mathbb{R}  = \re(-2xe^{-\imath\alpha}-2\bar{x}e^{\imath\alpha}) = -4\re(xe^{-\imath\alpha}).$$
Choose $\tilde{\alpha}$ such that
\begin{equation}\label{pivot-alpha}
e^{-\imath\tilde{\alpha}}=\textnormal{sgn}(\bar{x})=\frac{\bar{x}}{|x|}.
\end{equation}
Then
\begin{equation}\label{pivot-eta}
|\langle \grad f_{\Ham}(Z),Z\dot{R}(i,j,0,\tilde{\alpha}) \rangle_\mathbb{R}|=4\re(x\frac{\bar{x}}{|x|})=4|x|.
\end{equation}
Using relation~\eqref{pivot-gradfmax} we obtain
$$|\langle \grad f_{\Ham}(Z),Z\dot{R}(i,j,0,\tilde{\alpha}) \rangle_\mathbb{R}| \geq \frac{4}{\sqrt{4n^2-2n}}\|\grad f_{\Ham}(Z)\|_F=2\eta\|\grad f_{\Ham}(Z)\|_F.$$

If a unitary symplectic rotation $R(i,j,\phi,\alpha)$ is of the form~\eqref{srot3}, then
$$\dot{R}(i,j,0,\alpha)={\footnotesize\left[
  \begin{array}{ccccc|ccccc}
     &  &  &  &  &  &  &  &  &  \\
     &  &  &  &  &  &  &  & -e^{\imath \alpha} &  \\
     &  &  &  &  &  &  &  &  &  \\
     &  &  &  &  &  & -e^{-\imath \alpha} &  &  &  \\
     &  &  &  &  &  &  &  &  &  \\ \hline
     &  &  &  &  &  &  &  &  &  \\
     &  &  & e^{\imath \alpha} &  &  &  &  &  &  \\
     &  &  &  &  &  &  &  &  &  \\
     & e^{-\imath \alpha} &  &  &  &  &  &  &  &  \\
     &  &  &  &  &  &  &  &  &  \\
  \end{array}
\right]}$$
where in the matrix on the right-hand side all elements that are not explicitly given are zero.
The statement of the lemma follows same way as for the rotations~\eqref{srot2}.

If the rotation $R(i,j,\phi,\alpha)$ has the form~\eqref{srot1}, then
$$X^H\dot{R}(r,s,0,0)=
\left[
                         \begin{array}{ccccc}
                            &  &  &  &  \\
                            &  &  & -x &  \\
                            &  &  &  &  \\
                            & x &  &  &  \\
                            &  &  &  &  \\
                         \end{array}
                       \right]
\left[
                         \begin{array}{ccccc}
                            &  &  &  &  \\
                            &  &  & -1 &  \\
                            &  &  &  &  \\
                            & 1 &  &  &  \\
                            &  &  &  &  \\
                         \end{array}
                       \right],$$
where $x=x_{rs}$ is real.
Instead of~\eqref{pivot-eta} we get
$$|\langle \grad f(Z),Z\dot{R}(r,s,0,0) \rangle| = 2|x|.$$
Hence, from~\eqref{pivot-gradfmax} it follows
$$|\langle \grad f(Z),Z\dot{R}(r,s,0,\tilde{\alpha}) \rangle| \geq \frac{2}{\sqrt{4n^2-2n}}\|\grad f(Z)\|_F=\eta\|\grad f(Z)\|_F.$$
\end{proof}

\begin{Lemma}\label{lemma:cvgdelta-Ham}
Let $(Z_k,\ k\geq0)$ be the sequence generated by Algorithm~\ref{agm:Jacobi-Ham_closest}.
For every $\hat{Z}\in\U\cap\Symp$ with $\emph{grad} f_{\Ham}(\hat{Z})\neq0$,  there exist $\epsilon>0$ and $\delta>0$ such that
$$\|Z_k-\hat{Z}\|_F<\epsilon \quad \Rightarrow \quad f_{\Ham}(Z_{k+1})-f_{\Ham}(Z_k)\geq\delta.$$
\end{Lemma}

\begin{proof}
As $\grad f_{\Ham}(\hat{Z})\neq0$, there exists $\epsilon>0$ such that
\begin{equation}\label{cvgdelta-grad0}
\eta_1:=\min_{\|Z-\hat{Z}\|_F<\epsilon} \|\grad f_{\Ham}(Z)\|_F>0.
\end{equation}

For a fixed $k$ we define the differentiable function $h_k:\mathbb{R}\times\mathbb{R}\rightarrow\mathbb{R}$,
$$h_k(\phi,\alpha)=f_{\Ham}(Z_kR(i_k,j_k,\phi,\alpha)),$$
where $R(i_k,j_k,\phi,\alpha)$ is a unitary symplectic rotation as in Subsection~\ref{sec:rotations}.
As a part of Algorithm~\ref{agm:Jacobi-Ham_closest}, Algorithm~\ref{agm:angle} returns $\phi_k$ and $\alpha_k$ such that $f_\Ham$ is maximized.
Since $R(i_k,j_k,0,\alpha)=I$ for any $\alpha$, we have
\begin{equation}\label{cvgdelta-hk}
h_k(0,\alpha)=f_{\Ham}(Z_k) \quad \textnormal{and} \quad \max_{\phi,\alpha}h_k(\phi,\alpha)=h_k(\phi_k,\alpha_k)=f_{\Ham}(Z_{k+1}).
\end{equation}

Take $\tilde{\alpha}$ as in~\eqref{pivot-alpha} and define another function $H_k:\mathbb{R}\rightarrow\mathbb{R}$,
\begin{equation}\label{cvgdelta-Hdef}
H_k(\phi)=h_k(\phi,\tilde{\alpha}).
\end{equation}
The Taylor expansion of $H_k$ around $0$ yields
$$H_k(\phi_k) = H_k(0)+H_k'(0)\phi_k+\frac{1}{2}H_k''(\xi)\phi_k^2, \qquad 0<\xi<\phi_k.$$
Let $M=\max|H_k''(\xi)|<\infty$. Then
\begin{equation}\label{cvgdelta-hkphi}
H_k(\phi_k)-H_k(0) \geq H_k'(0)\phi_k-\frac{1}{2}M\phi_k^2.
\end{equation}

The derivative of $H_k$ is
$$H_k'(\phi) = \frac{\partial}{\partial\phi}f_{\Ham}(Z_kR(i_k,j_k,\phi,\tilde{\alpha}))
= \langle \grad f_{\Ham}(Z_kR(i_k,j_k,\phi,\tilde{\alpha})),Z_k\dot{R}(i_k,j_k,\phi,\tilde{\alpha}) \rangle_\mathbb{R},$$
and in particular,
$$H_k'(0) = \langle \grad f_{\Ham}(Z_k),Z_k\dot{R}(i_k,j_k,0,\tilde{\alpha}) \rangle_\mathbb{R}.$$
From Lemma~\ref{lemma:pivotcondition-Ham} and relation~\eqref{cvgdelta-grad0} we obtain
\begin{equation}\label{cvgdelta-hknul}
|H_k'(0)|\geq\eta\|\grad f_{\Ham}(Z_k)\| \geq \eta\min_{\|Z-\hat{Z}\|_F<\epsilon} \|\grad f_{\Ham}(Z)\|_F =\eta\eta_1.
\end{equation}
From~\eqref{cvgdelta-hk}, \eqref{cvgdelta-Hdef} and~\eqref{cvgdelta-hkphi}, for any $\phi$, we have
\begin{align}
f_{\Ham}(Z_{k+1})-f_{\Ham}(Z_k) & = h_k(\phi_k,\alpha_k)-h_k(0,\tilde{\alpha}) \geq h_k(\phi,\tilde{\alpha})-h_k(0,\tilde{\alpha}) \nonumber \\
& = H_k(\phi)-H_k(0) \geq H_k'(0)\phi-\frac{1}{2}M\phi^2. \label{cvgdelta-final}
\end{align}
Choose $\phi=\frac{H_k'(0)}{M}$. Finally, from,~\eqref{cvgdelta-final} and~\eqref{cvgdelta-hknul} we obtain
$$f_{\Ham}(Z_{k+1})-f_{\Ham}(Z_k) \geq \frac{H_k'(0)^2}{M}-\frac{H_k'(0)^2}{2M} = \frac{H_k'(0)^2}{2M}
\geq \frac{\eta^2\eta_1^2}{2M} = \delta.$$
\end{proof}

Now we can prove Theorem~\ref{tm:cvg-Ham} by contradiction.

\begin{proof}[Proof of Theorem~\ref{tm:cvg-Ham}]
Suppose that $\hat{Z}$ is an accumulation point of Algorithm~\ref{agm:Jacobi-Ham_closest}.
Then there is  a subsequence $\{Z_j\}, \ j\in K\subseteq\mathbb{N}$ such that $Z_j$ converges to $\hat{Z}.$

Assume that $\hat{Z}$ is not a stationary point of $f_{\Ham}$, that is $\grad f_{\Ham}(\hat{Z})\neq0$.
Then, for any $\epsilon > 0,$ there is $k_0\in K$ such that $\|Z_k-\hat{Z}\|<\epsilon$ for every $k>k_0$.
Lemma~\ref{lemma:cvgdelta-Ham} implies that $f_{\Ham}(Z_{k+1})-f_{\Ham}(Z_k)\geq\delta>0.$
Therefore, $f_{\Ham}(Z_k)\rightarrow\infty$ when $k\rightarrow\infty$.
Though, if $Z_k$ converges, $f_{\Ham}(Z_k)$ should converge, too. This gives a contradiction.
\end{proof}

\section{Numerical experiments}\label{sec:num}

We present some numerical experiments for the Hamiltonian and skew-Hamiltonian case.
We set up a random $2n\times2n$ Hamiltonian matrix $H$ as in~\eqref{eq:Hamsetup} by generating a random $n\times n$ matrix $H_{11}$ and random $n\times n$ Hermitian matrices $H_{12}$ and $H_{21}$.
Also, we set up a random $2n\times2n$ skew-Hamiltonian matrix $W$ as in~\eqref{eq:sHamsetup} using a random $n\times n$ matrix $W_{11}$ and random $n\times n$ skew-Hermitian matrices $W_{12}$ and $W_{21}$.
All tests were done in Matlab R2019b.

\begin{figure}[h!]
\begin{subfigure}{.45\textwidth}
    \centering
    \includegraphics[width=1\textwidth]{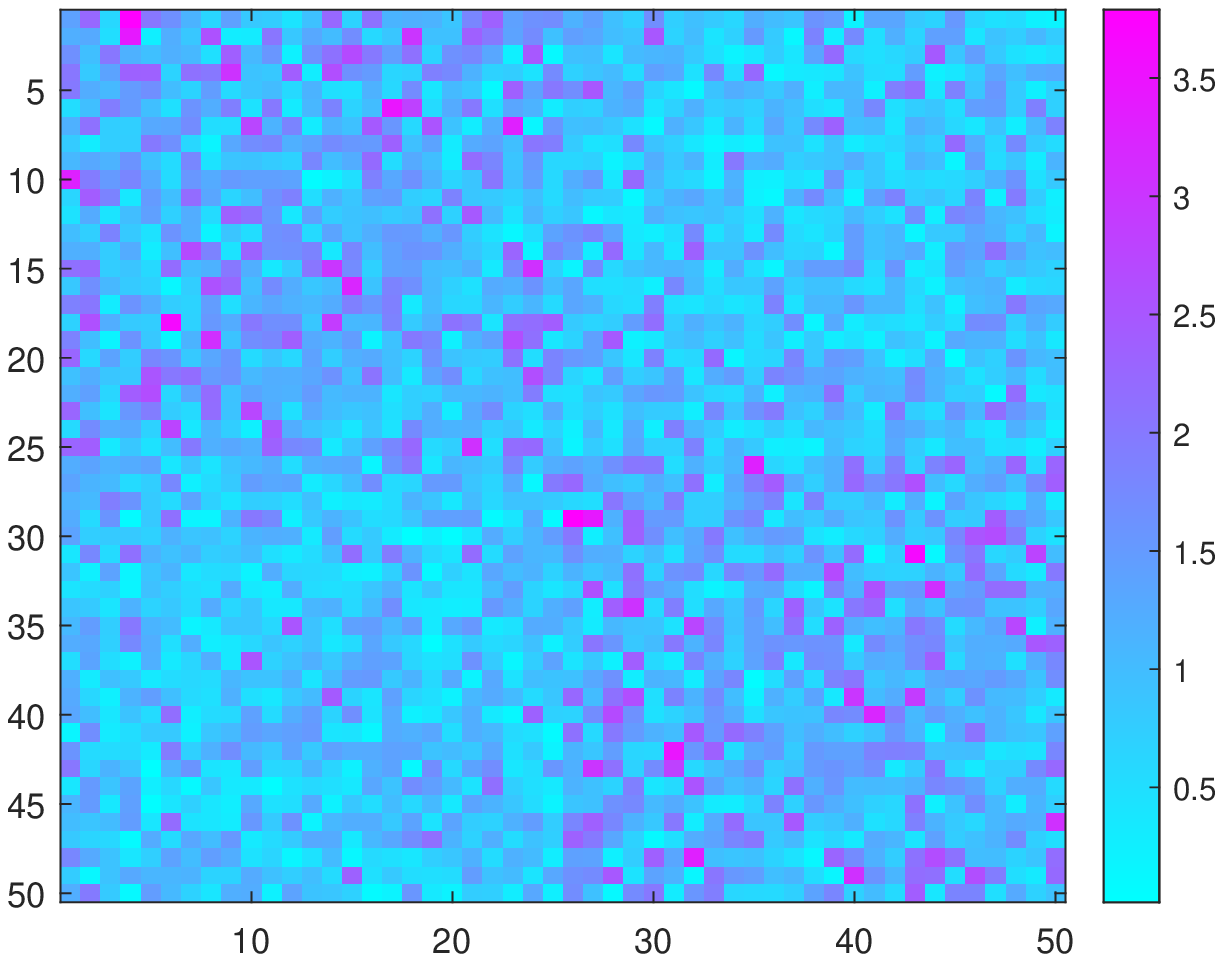}
    \caption{Hamiltonian matrix $H=H^{(0)}$}
\end{subfigure}
\begin{subfigure}{.45\textwidth}
    \centering
    \includegraphics[width=1\textwidth]{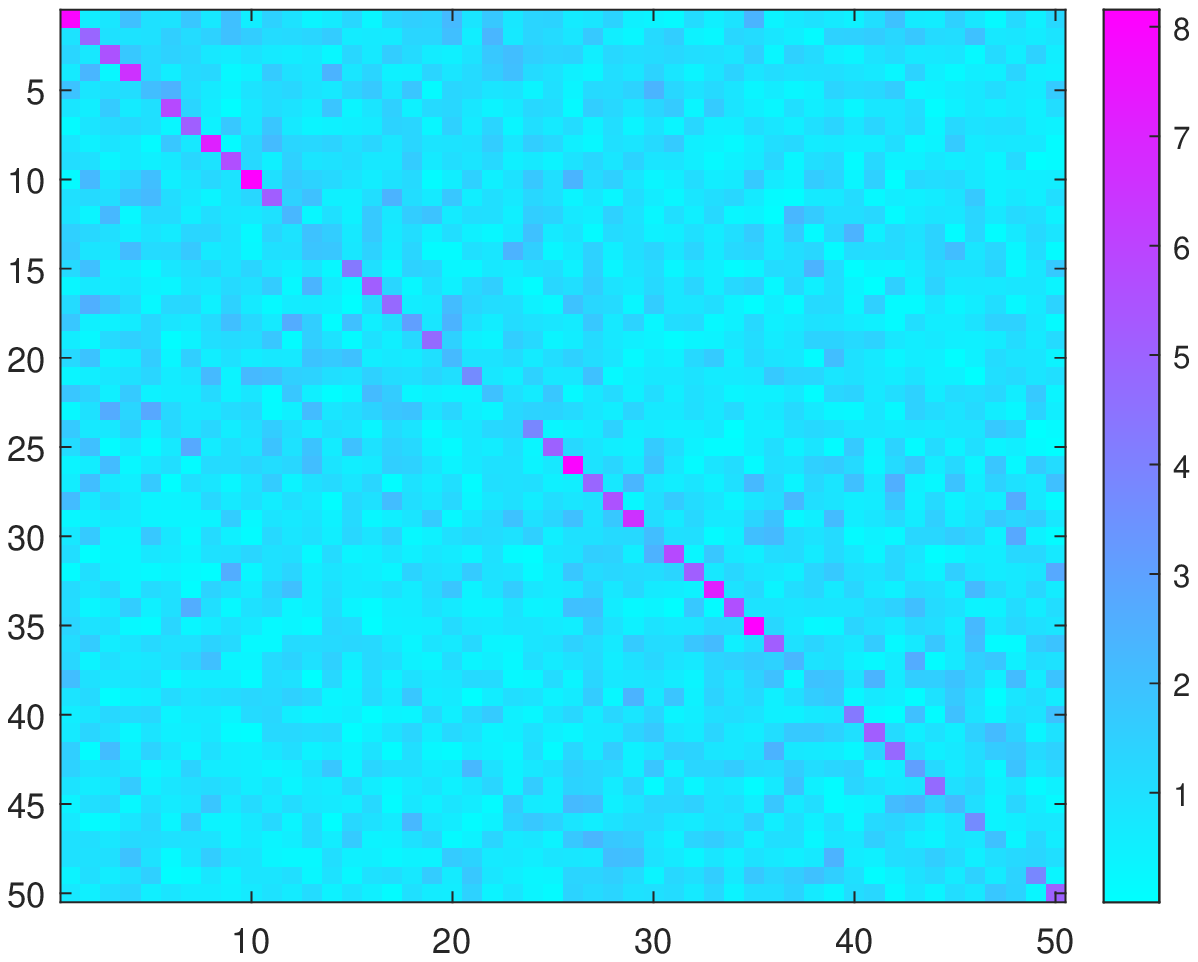}
    \caption{After $1$ iteration: $H^{(1)}$}
\end{subfigure}
\begin{subfigure}{.45\textwidth}
    \centering
    \includegraphics[width=1\textwidth]{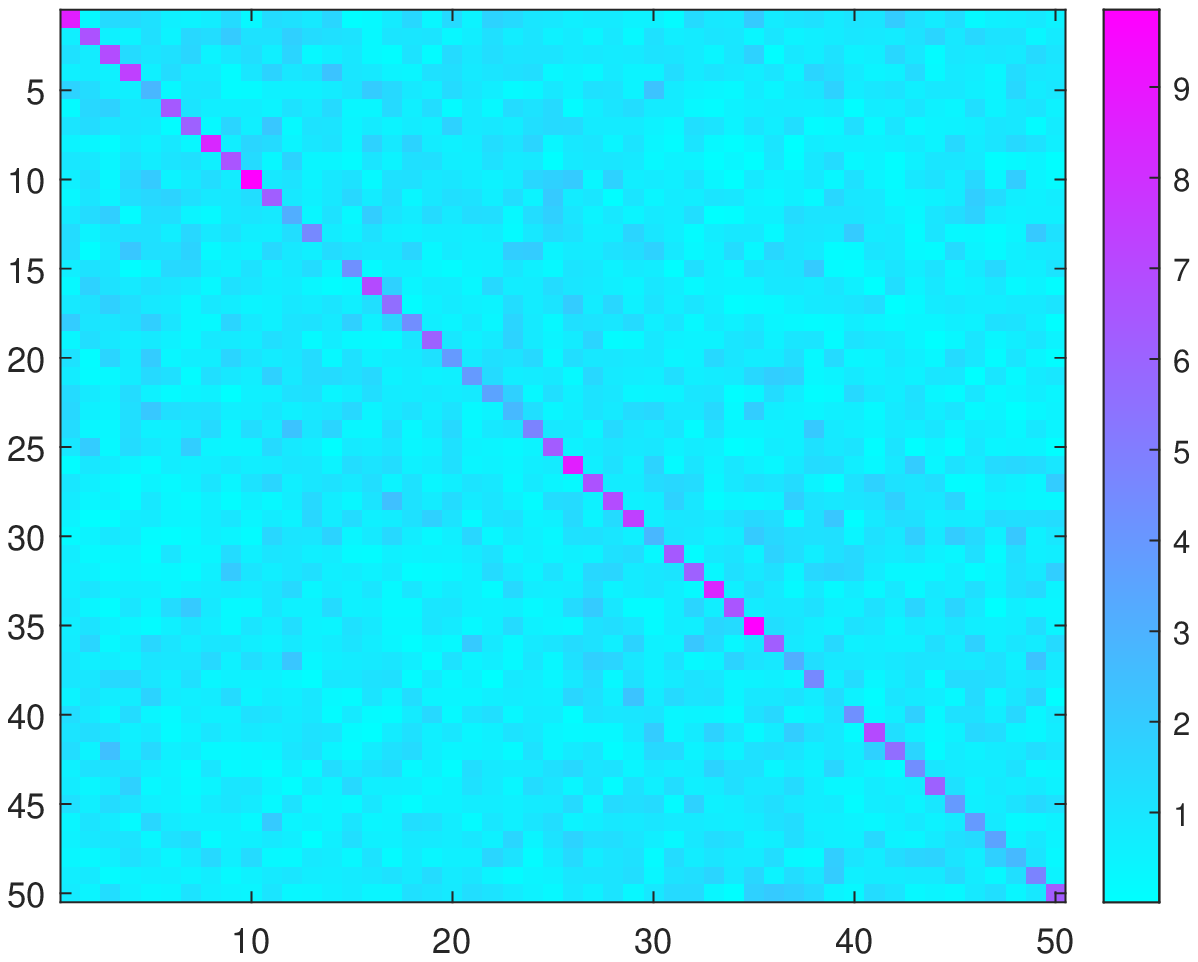}
    \caption{After $2$ iterations: $H^{(2)}$}
\end{subfigure}
\begin{subfigure}{.45\textwidth}
    \centering
    \includegraphics[width=1\textwidth]{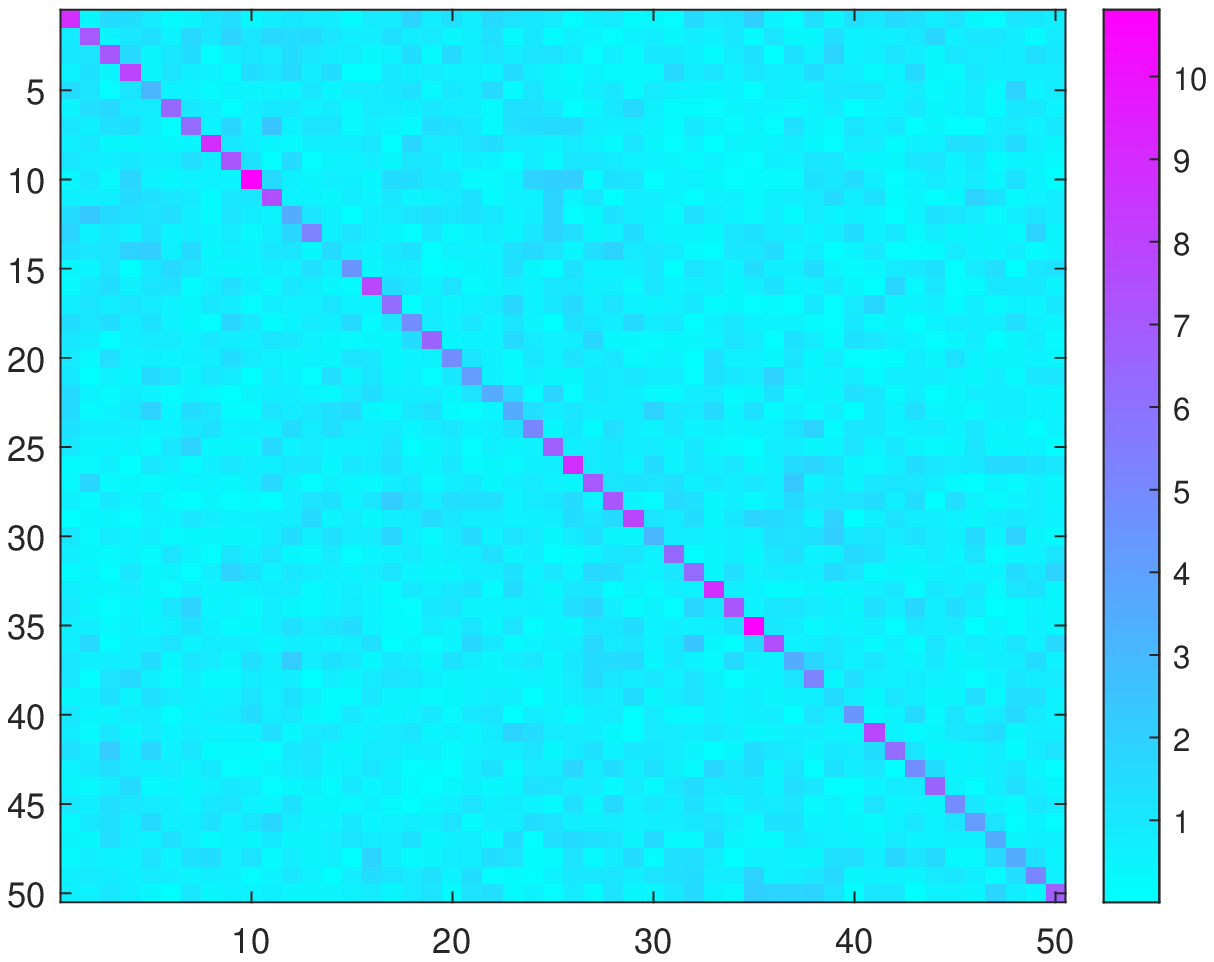}
    \caption{After $3$ iterations: $H^{(3)}$}
\end{subfigure}
    \caption{Change in the absolute value of matrix entries.}\label{fig:1}
\end{figure}

First, we see how the matrix norm moves to the diagonal during three iterations of Algorithm~\ref{agm:Jacobi-Ham_closest}.
In Figure~\ref{fig:1} the change in absolute value of the matrix entries is given. We start with a random $50\times50$ Hamiltonian matrix $H$ and show $H^{(k)}$, $k=1,2,3$. We can observe that the underlying matrix becomes diagonally dominant already after the first iteration. During the next iterations norm on the diagonal increases.

\begin{figure}[h!]
\begin{subfigure}{.45\textwidth}
    \centering
    \includegraphics[width=1\textwidth]{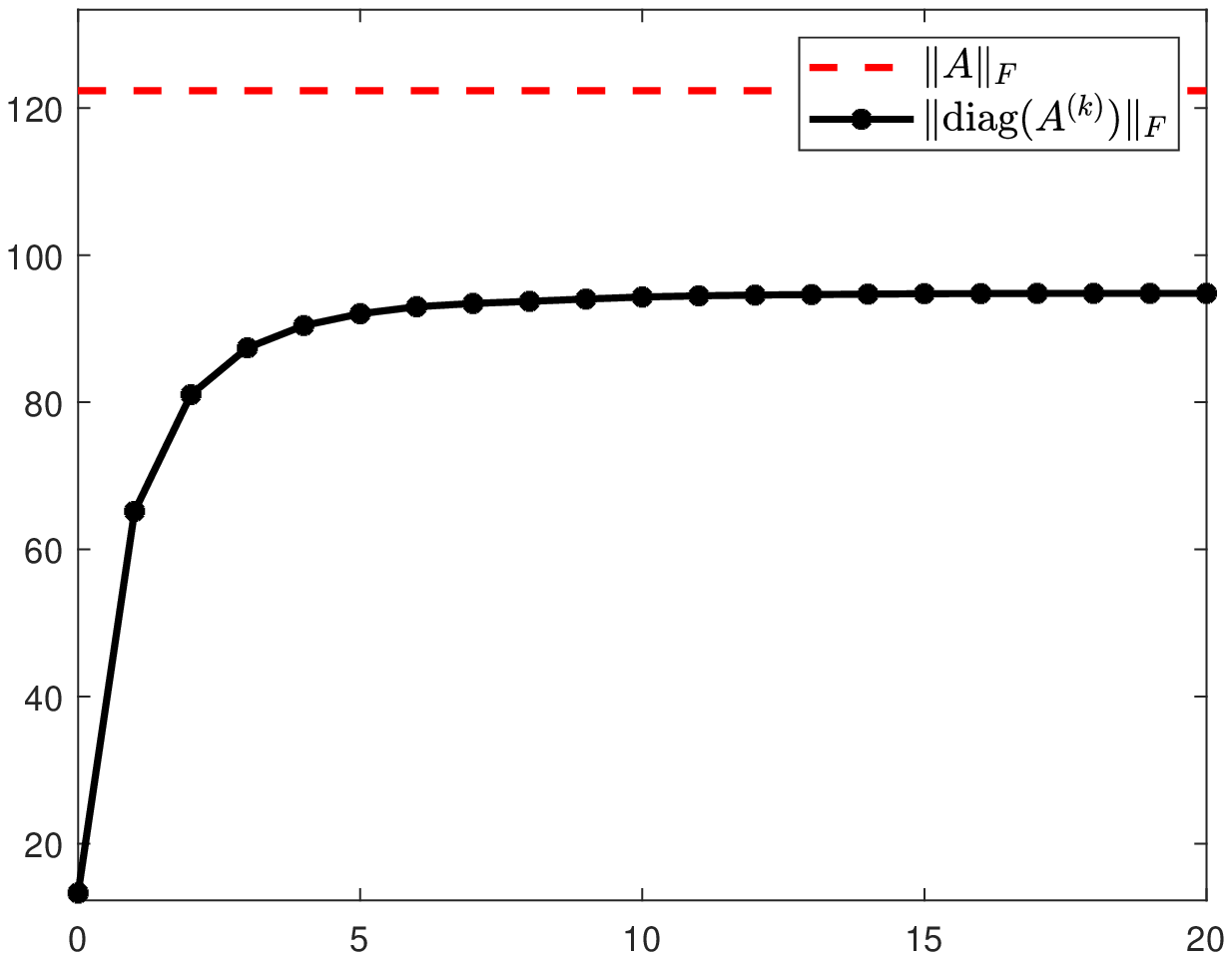}
    \caption{Matrix that can not be diagonalized using only symplectic rotations}
\end{subfigure}
\begin{subfigure}{.45\textwidth}
    \centering
    \includegraphics[width=1\textwidth]{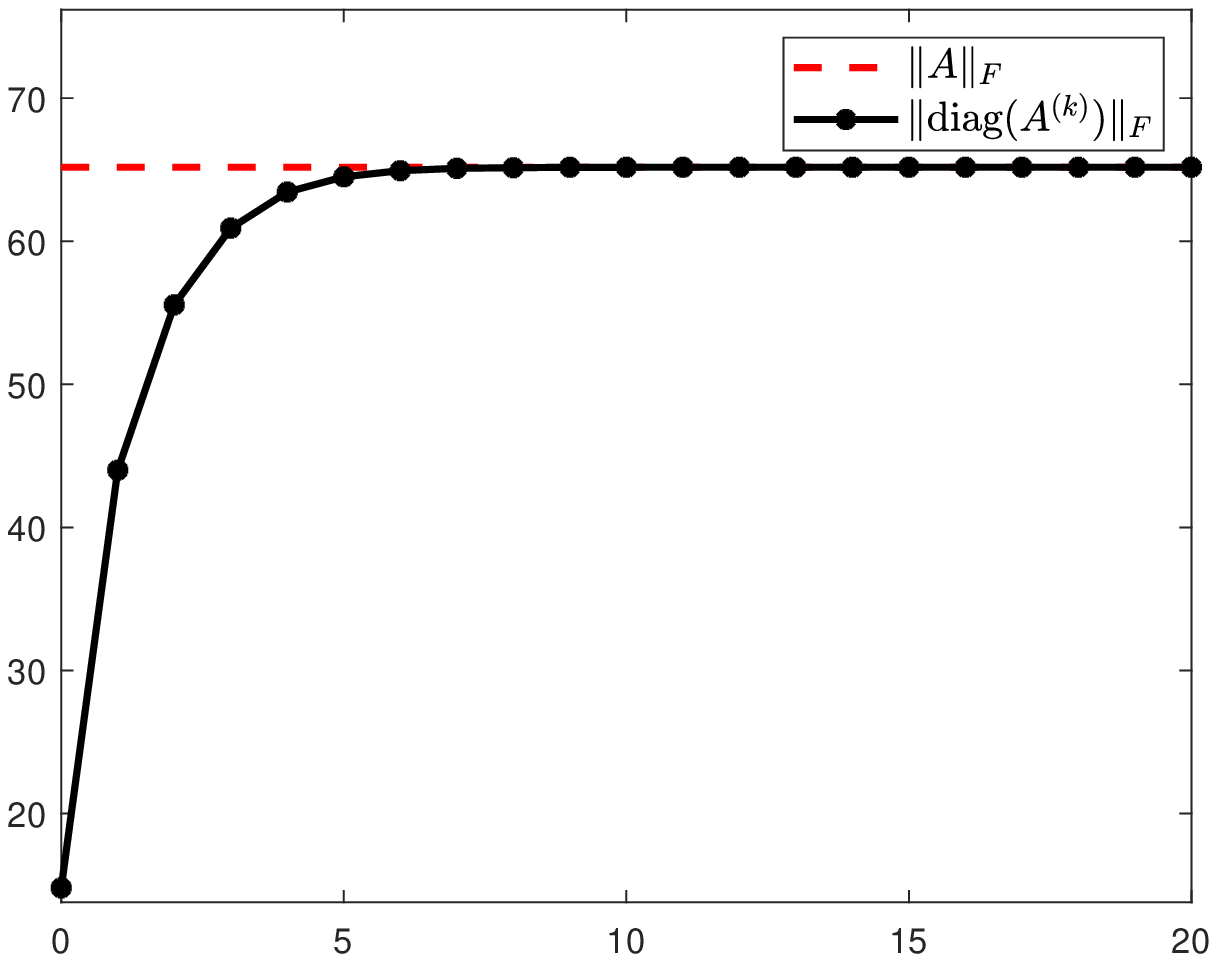}
    \caption{Matrix that can be diagonalized using only symplectic rotations}
\end{subfigure}
    \caption{Convergence of $\|\diag(A^{(k)})\|_F$, $A\in\Ham$.}\label{fig:2}
\end{figure}

In general, Hamiltonian matrix can not be diagonalized using symplectic rotations. Then Algorithm~\ref{agm:Jacobi-Ham_closest} will diagonalize it as much as possible. In Figure~\ref{fig:2} we show the convergence of $\|\diag(A^{(k)})\|_F$, $k=1,\ldots,20$, for two $100\times100$ Hamiltonian matrices, one that can not be diagonalized using only symplectic rotations and the other one that can. We see how $\|\diag(A^{(k)})\|_F$ approaches $\|A\|_F$. When $A$ can be diagonalized using only symplectic rotations, then complete norm of $A$ can be moved to its diagonal, so $\|\diag(A^{(k)})\|_F$ becomes equal to $\|A\|_F$.

\begin{figure}[h!]
\begin{subfigure}{.45\textwidth}
    \centering
    \includegraphics[width=1\textwidth]{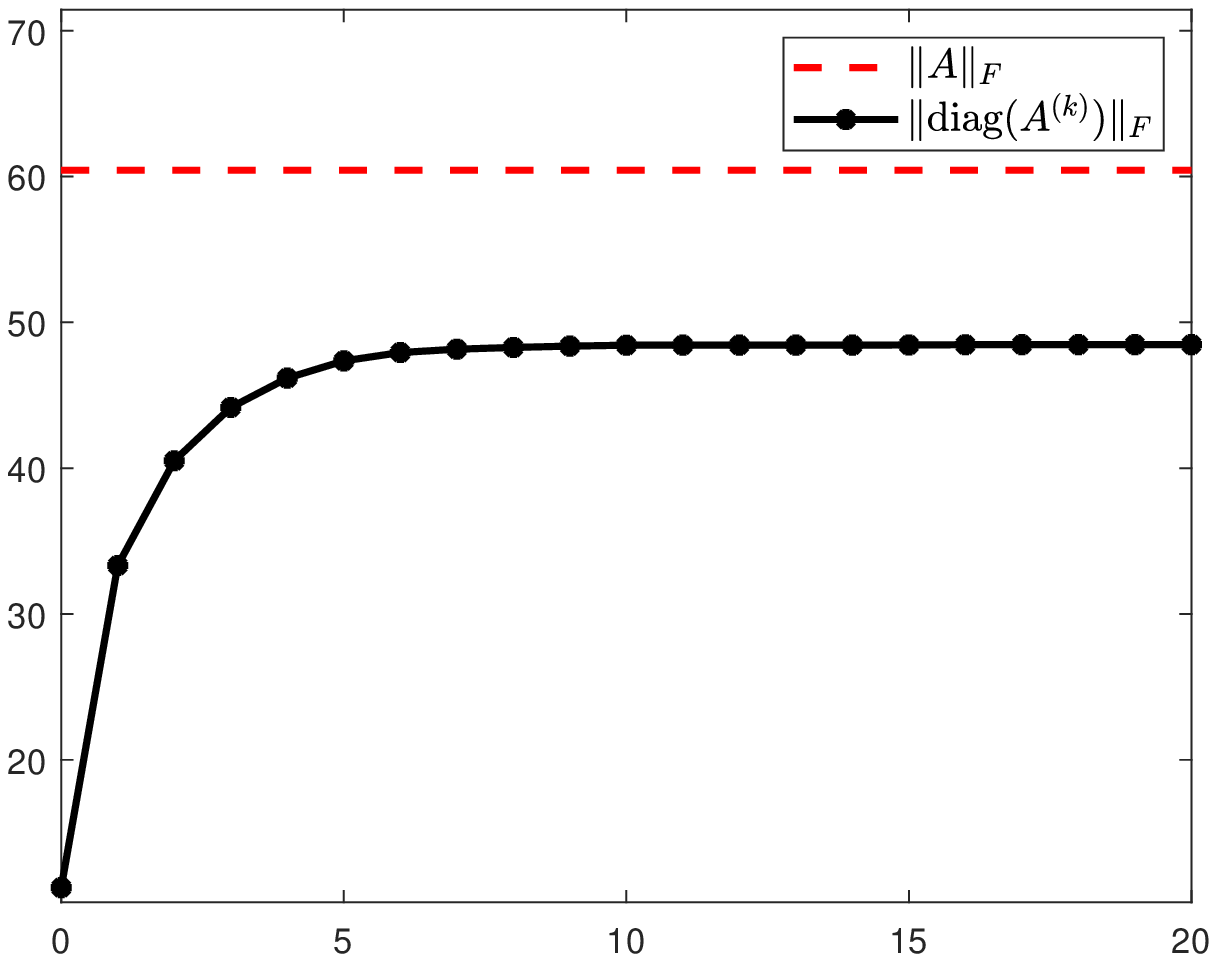}
    \caption{Matrix with no real eigenvalues}
\end{subfigure}
\begin{subfigure}{.45\textwidth}
    \centering
    \includegraphics[width=1\textwidth]{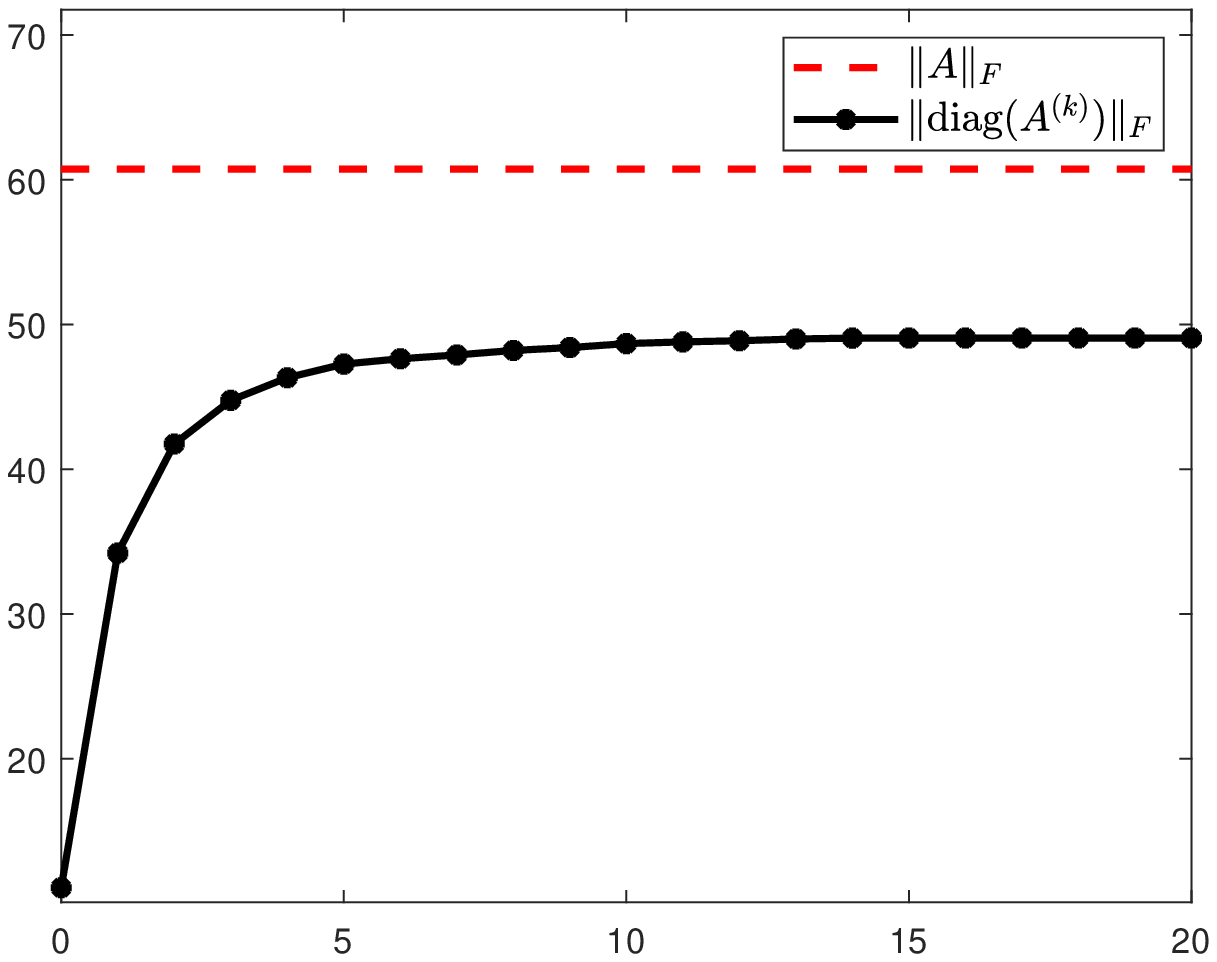}
    \caption{Matrix with some real eigenvalues}
\end{subfigure}
    \caption{Convergence of $\|\diag(A^{(k)})\|_F$, $A\in\sHam$.}\label{fig:3}
\end{figure}

Theorems~\ref{tm:maximizationH}, \ref{tm:maximizationsH}, \ref{tm:maximizationper} and~\ref{tm:maximizationpers} all include condition on the matrix eigenvalues. In Hamiltonian and perskew-Hermitian case, eigenvalues should not be purely imaginary, while in skew-Hamiltonian and per-Hermitian case they should not be real. Hamiltonian matrices from Figures~\ref{fig:1} and~\ref{fig:2} have all eigenvalues with non-zero real part, that is no purely imaginary eigenvalues.
Still, in practice, Algorithm~\ref{agm:Jacobi-Ham_closest} does not display any difference if the condition on the eigenvalues is not satisfied.
Figure~\ref{fig:3} gives the convergence of $\|\diag(A^{(k)})\|_F$ compared to $\|A\|_F$ for two $50\times50$ skew-Hamiltonian matrices, one with no and the other one with some real eigenvalues.

\begin{figure}[h!]
\begin{subfigure}{.45\textwidth}
    \centering
    \includegraphics[width=1\textwidth]{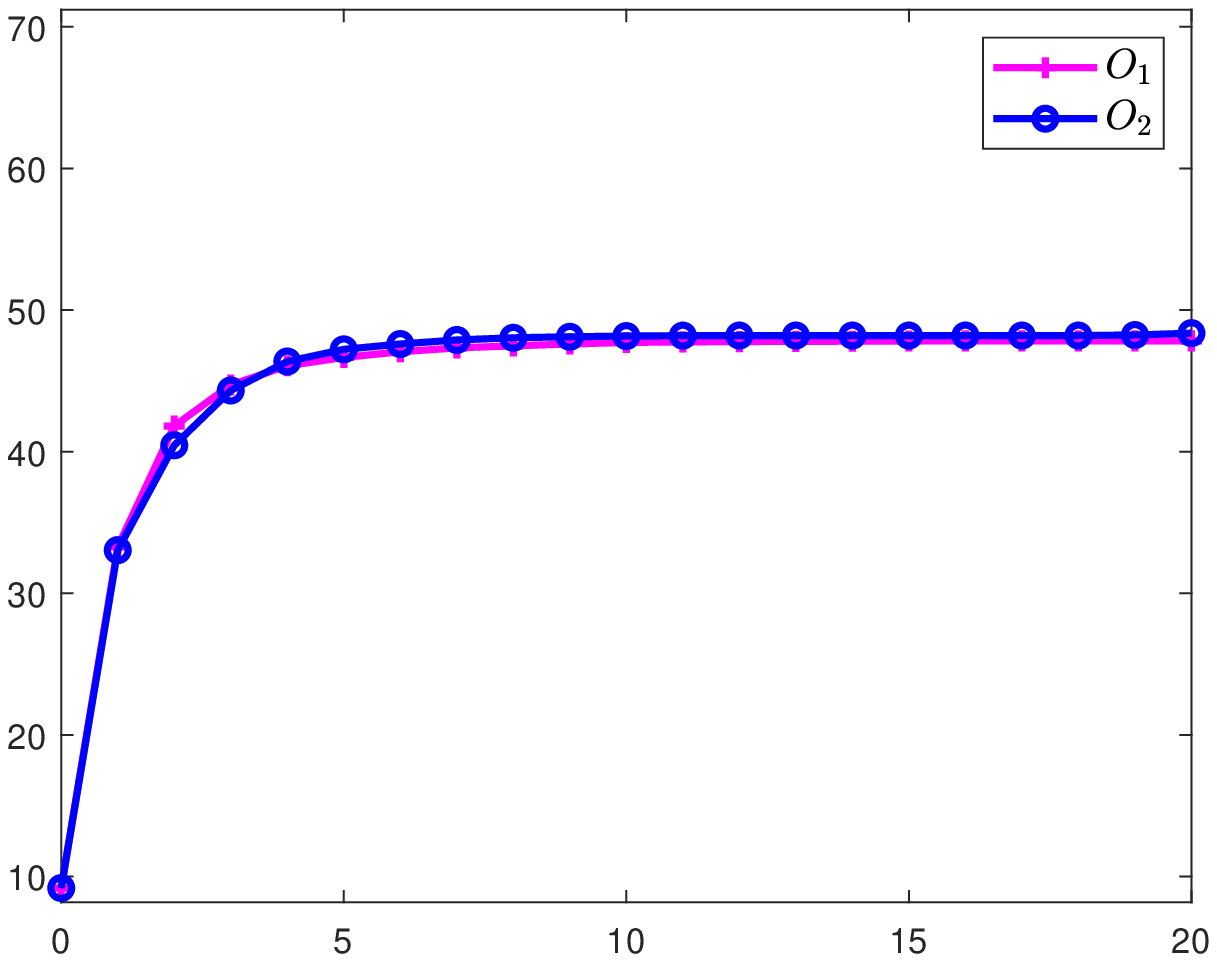}
    \caption{Matrix that can not be diagonalized by unitary symplectic transformations.}
\end{subfigure}
\begin{subfigure}{.45\textwidth}
    \centering
    \includegraphics[width=1\textwidth]{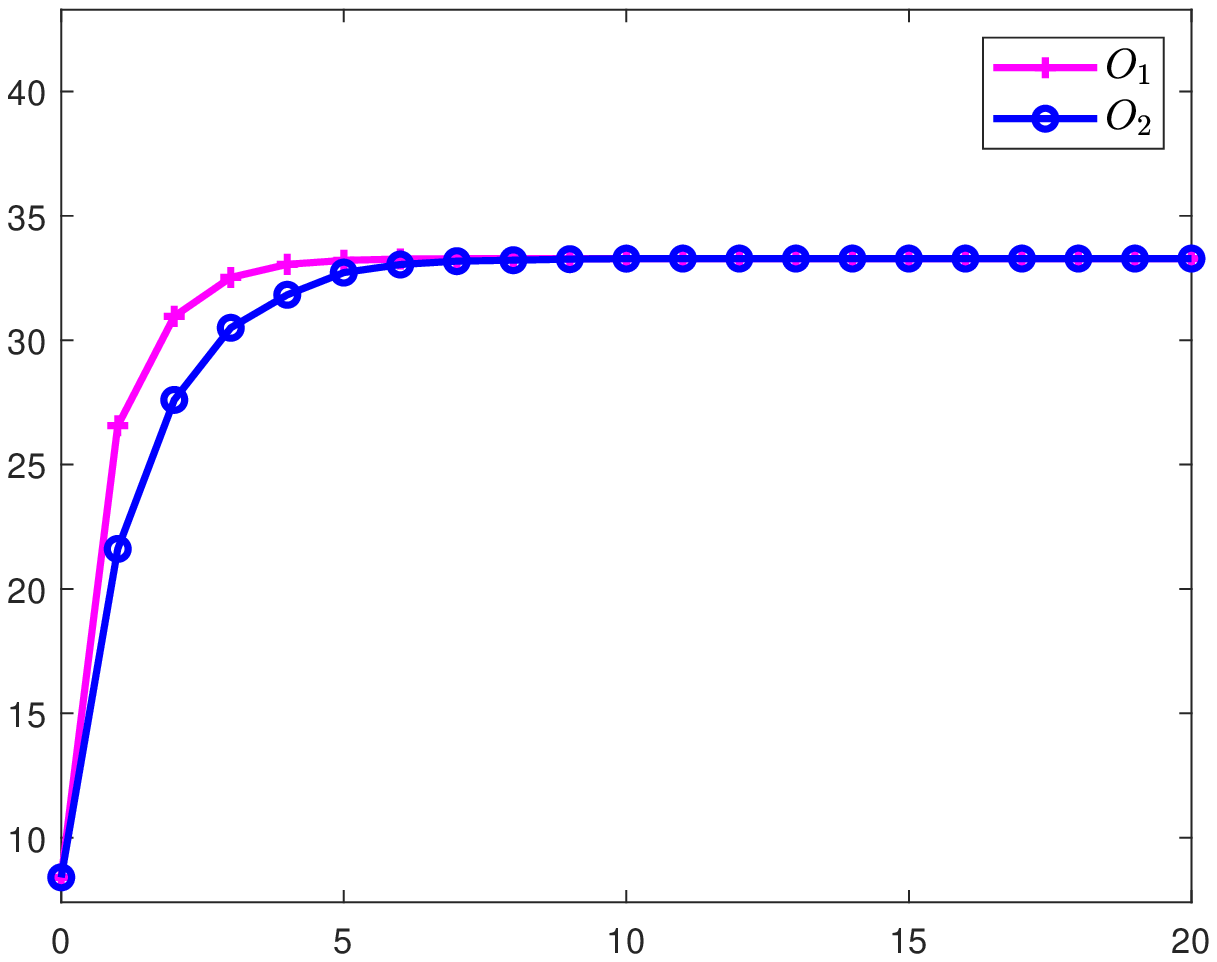}
    \caption{Matrix that can be diagonalized by unitary symplectic transformations.}
\end{subfigure}
    \caption{Convergence of $\|\diag(A^{(k)})\|_F$, $A\in\Ham$, for two different pivot strategies.}\label{fig:4}
\end{figure}

Recall that in Algorithm~\ref{agm:Jacobi-Ham_closest} any cyclic pivot ordering can be used since the convergence proof from Section~\ref{sec:cvgproof} does not depend on the ordering inside one sweep. In all previous examples we used pivot ordering
\begin{align*}
O_1 = & (1,2),(1,3),\ldots,(1,n),(2,3),\ldots,(2,n),\ldots,(n-1,n), \\
& (1,n+1),(2,n+2),\ldots,(n,2n), \\
& (1,n+2),(1,n+3),\ldots,(1,2n),(2,n+3),\ldots,(2,2n),\ldots,(n-1,2n).
\end{align*}
Now we will compare the convergence using two different cyclic orderings. The first one is $O_1$. The second one is ``bottom to top'' ordering from~\cite{Mehl08}. Keep in mind that we take pivot positions from the upper triangle, while in~\cite{Mehl08} they are taken from the lower triangle. This transforms ``bottom to top'' ordering into ``right to left'', meaning that instead of
$$(2n,1),(2n-1,1),\ldots,(2,1),(2n,2),\ldots,(3,2),\ldots,(2n,2n-1)$$
we have
\begin{equation}\label{ex:orderings}
(1,2n),(1,2n-1),\ldots,(1,2),(2,2n),\ldots,(2,3),\ldots,(2n-1,2n).
\end{equation}
Besides, because our algorithm uses double rotations, it does not take all pivot positions listed above, but its subset, as shown in~\eqref{pivotstrategy}. Thus, ``bottom to top'' ordering applied to our situation is a subset of~\eqref{ex:orderings} given by
\begin{align*}
O_2 = & (1,2n),(1,2n-1),\ldots,(1,n+1),(1,n),(1,n-1),\ldots,(1,2), \\
& (2,2n),(2,2n-1),\ldots,(2,n+2),(2,n),(2,n-1),\ldots,(2,3), \\
& (3,2n),\ldots,(3,n+3),(3,n),\ldots,(3,4),\ldots,(n-1,2n),(n-1,n),(n,2n).
\end{align*}
In Figure~\ref{fig:4} we present the convergence results for both orderings $O_1$ and $O_2$ on two random $50\times50$ Hamiltonian matrices, one that can not and one that can be completely diagonalized by unitary symplectic transformations.

\section*{Acknowledgements}

This work has been supported in part by Croatian Science Foundation under the project UIP-2019-04-5200 and by DAAD Short-term grant.
The author would like to thank Heike Fa{\ss}bender and Philip Saltenberger for useful discussions surrounding this problem.

\end{document}